\documentclass{amsart}
\usepackage{amssymb,amscd}
\advance\hoffset by -1in
\advance\voffset by -1in

\usepackage[cp1251]{inputenc}
\usepackage[russian]{babel}
\newtheorem{theorem}{Теорема}
\newtheorem{lemma}{Лемма}

\newtheorem{proposition}{Предложение}

\title{Теоремы о свободе для произведений групп}

\author{А.\,Ф.\,Красников}
\address{Омский государственный университет, г.\,Омск}
\email{phomsk@mail.ru}

\begin{document}

\maketitle

\section*{Введение}
\noindent Известная теорема о свободе Магнуса \cite{Mg} показывает, что если $R(x_1,\ldots,x_n)$ --- циклически несократимое слово в образующих
$x_1,\ldots,x_n$, содержащее $x_n$, то в группе $G=\langle x_1,\ldots,x_n;\text{ } R(x_1,\ldots,x_n)\rangle$ элементы
$x_1,\ldots,x_{n-1}$ являются свободными образующими порожденной ими подгруппы.\\
Подобное утверждение имеет место для групп с одним определяющим соотношением в многообразиях разрешимых и многообразиях нильпотентных групп данных ступеней \cite{Rm}, в многообразиях $\mathfrak{N}_c\mathfrak{A}$ \cite{Ya} (здесь $\mathfrak{A}$ --- многообразие абелевых групп, $\mathfrak{N}_c$ --- многообразие нильпотентных групп ступени, не превосходящей $c$), в полинильпотентных многообразиях  \cite{Km}.\\
Пусть $F=\underset{i\in I}\ast A_i$ --- свободное
произведение нетривиальных групп $A_i~(i\in I)$. Эндоморфизм $\varphi$ группы $F$ будем называть элементарным, если найдется $K_\varphi\subseteq I$
такое, что $\mbox{если }a\in A_i~\mbox{то }\varphi(a)=1~(i\in K_\varphi),~\varphi(a)=a~(i\not\in K_\varphi)$.
В настоящей работе доказана теорема:
\begin{theorem}\label{tm1_pr_gr}
Пусть $F$ --- свободное
произведение нетривиальных групп $A_i~(i=1,\ldots,n)$ $(n>2)$,
$1\neq N$ --- нормальная, допускающая элементарные эндоморфизмы группы $F$ подгруппа в $F$ такая, что $N\cap A_i=1~(i=1,\ldots,n)$, $F/N$ --- упорядочиваемая группа. Пусть, далее,
\begin{eqnarray}\label{end_gr}
N=N_{11} > \ldots > N_{1,m_1+1}=N_{21} > \ldots > N_{s,m_s+1},
\end{eqnarray}
где $N_{ij}$ --- j-й член нижнего центрального ряда группы $N_{i1}$, $s \geqslant 1$,
$H$ --- свободное произведение групп $A_i~(i=1,\ldots,n-1)$, $r\in N_{1t}\backslash N_{1,t+1}$, $R$ --- нормальная подгруппа, порожденная в группе $F$ элементом $r$.\\
Если (и только если) элемент $r$ не сопряжен по модулю $N_{1,t+1}$ ни с каким элементом из $H$, то
$H\cap RN_{kl} = H\cap N_{kl}$, где $N_{kl}$ --- произвольный член ряда {\rm (\ref{end_gr})}.
\end{theorem}
\noindent Теорема Романовского \cite{Rm1} показывает, что если группа $G$ задана с помощью образующих
элементов $x_1,\ldots,x_n$ и определяющих соотношений $r_1,\ldots,r_m$, где $n>m$, то среди
данных образующих существуют такие $n-m$ элементов $x_{i_1},\ldots,x_{i_{n-m}}$, которые свободно порождают свободную подгруппу в $G$.\\
Подобное утверждение имеет место для разрешимых групп, т.е. для групп, заданных с помощью образующих
элементов и определяющих соотношений в многообразиях разрешимых групп ступени разрешимости не выше $k$ \cite{Rm1},
для про-$p$-групп, разрешимых про-$p$-групп, нильпотентных про-$p$-групп и нильпотентных абстрактных групп \cite{Rm2},
а также для групп, заданных образующими и соотношениями в некоторых других многообразиях абстрактных групп \cite{Rm3,Rm4,Rm5}
и для свободных и разрешимых произведений групп с соотношениями \cite{Rm6}.\\
В настоящей работе доказана теорема:
\begin{theorem}\label{tm2_pr_gr}
Пусть $F$ --- свободное
произведение нетривиальных групп $A_i~(i=1,\ldots,n)$,
$1\neq N$ --- нормальная, допускающая элементарные эндоморфизмы группы $F$ подгруппа в $F$ такая, что $N\cap A_i=1~(i=1,\ldots,n)$, $F/N$ --- правоупорядочиваемая разрешимая группа. Пусть, далее,
\begin{eqnarray}\label{end_gr_1}
N=N_{11} > \ldots > N_{1,m_1+1}=N_{21} > \ldots > N_{s,m_s+1},
\end{eqnarray}
где $N_{ij}$ --- j-й член нижнего центрального ряда группы $N_{i1}$, $s \geqslant 1$,
$r_1,\ldots,r_m$ --- элементы из $N$ $(m<n)$, $R$ --- нормальная подгруппа, порожденная в группе $F$ элементами $r_1,\ldots,r_m$.
Тогда найдется $J \subseteq \{1,2,\ldots,n\}$, $|J|\geqslant n-m$,
такое, что для подгруппы $H$ группы $F$, порожденной группами $A_i,~i\in J$, справедлива формула
$H\cap RN_{kl} = H\cap N_{kl}$, где $N_{kl}$ --- произвольный член ряда {\rm (\ref{end_gr_1})}.
\end{theorem}

\section{Некоторые свойства производных Фокса}
\noindent Пусть $F=(\underset{i\in I}\ast A_i)\ast G$ --- свободное
произведение нетривиальных групп $A_i~(i\in I)$ и свободной группы $G$ с базой
$\{g_j | j\in J\}$; $N$ --- нормальная подгруппа в $F$.\\
Обозначим через ${\bf Z}[F]$ целочисленное групповое кольцо группы
$F$. Дифференцированием кольца ${\bf Z}[F]$ называется отображение
$\partial:{\bf Z}[F]\to {\bf Z}[F]$ удовлетворяющее условиям
\begin{gather}
\partial(u+v)=\partial(u)+\partial(v),\notag\\
\partial(uv)=\partial(u)v+ \varepsilon (u)\partial(v)\notag
\end{gather}
для любых $u,~v\in {\bf Z}[F]$, где $\varepsilon$ --- гомоморфизм
тривиализации $F\to 1$, продолженный по линейности на ${\bf Z}[F]$.\\
Обозначим через $D_k~(k\in I\cup J)$ производные Фокса кольца
${\bf Z}[F]$ --- дифференцирования, однозначно определяемые
условиями:
\begin{gather}
D_j(g_j)=1~(j\in J),~D_k(g_j)=0,~k\neq j;\notag\\
\mbox{если}~a_i\in A_i~(i\in
I),~\mbox{то}~D_i(a_i)=a_i-1,~D_k(a_i)=0,~k\neq i.\notag
\end{gather}
Пусть $\Delta_i$ --- фундаментальный идеал кольца ${\bf Z}A_i$. Имеют место формулы:
\begin{gather}
\mbox{если}~u_i\in {\bf Z}[F]\cdot\Delta_i~(i\in
I)~\mbox{то}~D_i(u_i)=u_i,~D_k(u_i)=0,~k\neq i;\notag\\
D_j((g_j-1)u)=u,~D_k((g_j-1)u)=0,~k\neq j;\notag\\
D_k(f^{-1}nf)\equiv D_k(n)f\mod{{\bf Z}[F]\cdot (N-1)};~D_k(f^{-1})=-D_k(f)f^{-1};\notag\\
u-\varepsilon(u)=\sum_{i\in I} D_i(u)+\sum_{j\in J}
(g_j-1)D_j(u);\label{1tm01}
\end{gather}
где $u\in {\bf Z}[F],~f\in F,~n\in N$.\\
Пусть $G$ --- группа; $A$, $B$ --- подмножества множества элементов группы $G$; $C$, $D$ --- подмножества множества элементов кольца ${\bf Z}[G]$. Обозначим через $\textup{гр}(A)$ подгруппу, порожденную $A$ в $G$, через $A^G$ --- нормальную подгруппу,
порожденную $A$ в $G$, через $\gamma_k G$ --- $k$-й член нижнего центрального ряда группы $G$. Если $x,~y$ --- элементы $G$, то положим $[x,y]=x^{-1}y^{-1}xy$, $x^y=y^{-1}xy$.
Через $AB$ обозначим множество произведений вида $ab$, где $a,~b$ пробегают соответственно элементы $A$, $B$, через $[A,B]$ --- подгруппу группы $G$, порожденную всеми $[a,b]$, $a\in A$, $b\in B$. Через $CD$ обозначим множество сумм произведений вида $cd$, где $c,~d$ пробегают соответственно элементы $C$, $D$.

\begin{theorem}\label{tm1_gr}
Пусть $F$ --- свободное произведение нетривиальных групп
$A_i~(i\in I)$ и свободной группы с базой $\{g_j | j\in J\}$,
$N$--- нормальная подгруппа в $F$, $K\subseteq I\cup J$, $F_K$
 --- подгруппа в $F$, порожденная $\{g_j| j\in K\cap J\}$ и
$\{A_i| i\in K\cap I\}$; $v\in F$, $M = \textup{гр
}((N\cap A_i)^{F}|i\in I)[N,N]$. Тогда
\begin{eqnarray}\label{1tm02}
D_k(v)\equiv ~0\mod{{\bf Z}[F]\cdot (N-1)},~k\in (I\cup J)\setminus K,
\end{eqnarray}
если и только если найдется элемент $\hat{v}\in F_K$ такой, что $v\hat{v}^{-1}\in (F_K\cap N)^FM$.
\end{theorem}
\begin{proof}
Непосредственно проверяется, что из $v\hat{v}^{-1}\in (F_K\cap N)^FM$, $\hat{v}\in F_K$, следуют сравнения (\ref{1tm02}).
Необходимо доказать обратное.\\
Через $\mathfrak{N}$ будем обозначать идеал, порожденный в ${\bf Z}[F]$ элементами
$\{n-1| n\in N\}$.
Формулы (\ref{1tm01}), (\ref{1tm02}) показывают, что
\begin{eqnarray}\label{1tm01-02}
{v-1}\equiv     \sum_{i\in K\cap I} D_i(v) + \sum_{j\in K\cap
J}(g_j-1)D_j(v)\mod{\mathfrak{N}}.
\end{eqnarray}
Обозначим через $\Delta$ фундаментальный идеал кольца ${\bf Z}[F_K]$.
Из (\ref{1tm01-02}) следует, что при естественном гомоморфизме
${\bf Z}[F]\to {\bf Z}[F/N]$ образ элемента $v-1$ принадлежит
образу $\Delta\cdot {\bf Z}[F]$.
Следовательно, элементы $1$ и $v$ принадлежат одному правому классу смежности группы
$F$ по подгруппе $F_KN$, т.е. найдется $\hat{v}\in F_K$ такой, что $v\hat{v}^{-1}\in N$.
Элемент $v\hat{v}^{-1}$ обозначим через $w$.
Предположим, что $w\notin (F_K\cap N)^FM$ и приведем это предположение к противоречию.
Пусть $X$ --- множество элементов из $F$, полученное объединением $\{g_j
| j\in J\}$ с множеством элементов групп $A_i~(i\in I)$; $u\to
\bar u$ --- функция, выбирающая правые шрайеровы представители $F$ по $N$,
$S$ --- множество выбранных представителей. Тогда $N = \text{гр }(sx{\overline{sx}}^{-1}|s\in S,~x\in X)$ (доказательство
см., например, в \cite{KrMr}). Следовательно, найдутся $s_ix_i{\overline{s_ix_i}}^{-1},~s_i\in S,~x_i\in
X~(i=1,\ldots,l),~k_1,\ldots,k_l$ --- целые, отличные от нуля
числа такие, что $w\equiv (s_1x_1{\overline{s_1x_1}}^{-1})^{k_1}\ldots (s_lx_l{\overline{s_lx_l}}^{-1})^{k_l}\mod{(F_K\cap N)^FM}$ и элемент $w$ нельзя представить по модулю $(F_K\cap N)^FM$ в виде произведения степеней меньшего чем $l$ числа элементов вида
$sx{\overline{sx}}^{-1},~s\in S,~x\in X$. Через $w_i$ будем обозначать элемент $s_ix_i{\overline{s_ix_i}}^{-1}$,
i=1,\ldots,l.
Отметим, что
\begin{eqnarray}\label{1tm03}
\sum_{i=1}^l k_i D_q(w_i)\equiv 0 ~\mod{\mathfrak{N}},~q\in (I\cup J)\setminus K.
\end{eqnarray}
Предположим, что в $\{x_1,\ldots,x_l\}$ есть не принадлежащий $F_K$ элемент. Пусть это будет $x_1$.
Приведем это предположение к противоречию.
Выберем $q\in (I\cup J)\setminus K$ такое, что либо $x_1=g_q$ либо $x_1\in A_q$.
Рассмотрим случай $x_1=g_q$. Если $x_i=g_q$, то
\begin{eqnarray}\label{1tm04}
D_q(w_i)\equiv D_q(s_i){s_i}^{-1}+
{\overline{s_ix_i}}^{-1}-D_q(\overline{s_ix_i}){\overline{s_ix_i}}^{-1}~\mod{\mathfrak{N}},
\end{eqnarray}
если $x_i\neq g_q$, то
\begin{eqnarray}\label{1tm05}
D_q(w_i)\equiv
D_q(s_i){s_i}^{-1}-D_q(\overline{s_ix_i}){\overline{s_ix_i}}^{-1}~\mod{\mathfrak{N}}.
\end{eqnarray}
Пусть $t\in S$ и $t=u{g_q}^\varepsilon u_1$, где $\varepsilon = \pm 1$ и слово $u_1$ не
содержит в своей записи буквы $g_q$. Тогда
\begin{eqnarray*}
D_q(t){t}^{-1}=
D_q(u){u}^{-1}+D_q({g_q}^\varepsilon)({u{g_q}^\varepsilon})^{-1},
\end{eqnarray*}
т.е. $D_q(t){t}^{-1}$ --- сумма элементов вида $D_q({g_q}^\varepsilon)({u{g_q}^\varepsilon})^{-1}$, $u{g_q}^\varepsilon\in S$.
Покажем, что
\begin{eqnarray*}
{\overline{s_1x_1}}^{-1}\equiv \pm
D_q({g_q}^\varepsilon)({u{g_q}^\varepsilon})^{-1}~\mod~N.
\end{eqnarray*}
Предположим противное. При $\varepsilon = 1$ будем иметь ${\overline{s_1x_1}}^{-1}\equiv (u{g_q})^{-1}~\mod~N$.
Тогда $s_1=u$ и
$\overline{s_1x_1}=s_1x_1$. Пришли к
противоречию.\\
При $\varepsilon = -1$ будем иметь ${\overline{s_1x_1}}^{-1}\equiv {u}^{-1}~\mod~N$.\\
Тогда $s_1=ux_1^{-1}$ и
$s_1x_1{\overline{s_1x_1}}^{-1}=1$.
Пришли к противоречию.\\
Теперь можно утверждать, что из (\ref{1tm03}), (\ref{1tm04}),
(\ref{1tm05}) следует существование $i$ такого, что $i\neq
1,~x_i=g_q$ и
${\overline{s_ix_i}}^{-1}={\overline{s_1x_1}}^{-1}$. Тогда $s_i=s_1$ и потому $w_i=w_1$. Снова пришли к противоречию.\\
Рассмотрим оставшийся случай $x_1\in A_q$. Обозначим через $L$ подмножество в
$\{1,\ldots,l\}$, состоящее из индексов тех элементов $x_i$, которые принадлежат
$A_q$. Выберем элементы  $u_i$, $v_i$, $a_i$, $i\in L$, так, чтобы было:
$u_i$, $v_i$ --- элементы из $S$ и не кончаются символом из $A_q$, $a_i\in A_q$,
$w_i=u_ia_i{v_i}^{-1}$, $i\in L$.\\
Если $i\in L$, то
\begin{eqnarray}\label{1tm06}
D_q(w_i)\equiv D_q(u_i){u_i}^{-1}+
{u_i}^{-1}-{v_i}^{-1}-D_q(v_i){v_i}^{-1}~\mod{\mathfrak{N}}.
\end{eqnarray}
Если $i\notin L$, то
\begin{eqnarray}\label{1tm07}
D_q(w_i)\equiv
D_q(s_i){s_i}^{-1}-D_q(\overline{s_ix_i}){\overline{s_ix_i}}^{-1}~\mod{\mathfrak{N}}.
\end{eqnarray}
Предположим, что $\sum_{j\in L} k_j({u_j}^{-1}-{v_j}^{-1})\neq 0$.
Пусть $y\in F$ и $y=y_1y_2$, где слово $y_2$ не содержит в своей записи
символов из $A_q$. Тогда $D_q(y){y}^{-1}=D_q(y_1){y_1}^{-1}$,
поэтому из (\ref{1tm03}), (\ref{1tm06}), (\ref{1tm07}) следует, что найдутся кончающиеся символом из $A_q$ попарно
различные элементы $z_1,\ldots,z_n$ из $S$ и целые, не равные нулю,
числа $m_1,\ldots,m_n$ такие, что
\begin{eqnarray}\label{1tm08}
\sum_{j\in L} k_j({u_j}^{-1}-{v_j}^{-1})+\sum_{i=1}^n m_i
D_q(z_i){z_i}^{-1}\equiv 0 ~\mod{\mathfrak{N}}.
\end{eqnarray}
Пусть $z_i=\hat{z}_i a_i$, где $a_i\in A_q$. Тогда
\begin{eqnarray}\label{1tm09}
D_q(z_i){z_i}^{-1} = D_q(\hat{z}_i){\hat{z}_i}^{-1} +
{\hat{z}_i}^{-1} - {z_i}^{-1}.
\end{eqnarray}
Из (\ref{1tm08}), (\ref{1tm09}) вытекает $m_i=0$ для тех $i$,
для которых слова $z_i$ имеют в своей записи максимальное число
вхождений символов из $A_q$. Пришли к противоречию.\\
Нами доказано, что $\sum_{j\in L} k_j({u_j}^{-1}-{v_j}^{-1})= 0$.
Так как $w_j\notin M$ и $u_ja_j{v_j}^{-1}\in N$, то
${u_j}^{-1}\neq {v_j}^{-1}~(j\in L)$.  Поэтому для любого элемента
$b\in \{u_i,~v_i\}$ найдется $\{u_j,~v_j\}~(j\neq i,~i,j\in L)$
такое, что
$b\in \{u_j,~v_j\}$.\\
Следовательно, найдутся $\{u_{i_1},~v_{i_1}\},
\ldots,\{u_{i_k},~v_{i_k}\}~(i_1,\ldots,i_k$ - попарно различные
элементы из $L$), попарно различные элементы $b_1,\ldots,b_k$ из
$S$ и $\varepsilon_i=\pm 1~(i=1,\ldots,k)$ такие, что
\begin{gather}
\{u_{i_k},~v_{i_k}\}\cap\{u_{i_1},~v_{i_1}\}=b_1,~\{u_{i_{j-1}},~v_{i_{j-1}}\}\cap\{u_{i_j},~v_{i_j}\}=b_j,~j=2,\ldots,k,\notag\\
{w_{i_1}}^{\varepsilon_1}\ldots
{w_{i_k}}^{\varepsilon_k}=b_1{x_{i_1}}^{\varepsilon_1}\ldots
{x_{i_k}}^{\varepsilon_k}{b_1}^{-1}.\notag
\end{gather}
Из ${x_{i_1}}^{\varepsilon_1}\ldots {x_{i_k}}^{\varepsilon_k}\in N\cap A_q$
следует, что
${w_{i_1}}^{\varepsilon_1}\ldots {w_{i_k}}^{\varepsilon_k}\in M$ в противоречии с выбором $w_1,\ldots,w_l$.\\
Полученные противоречия показывают, что $\{x_1,\ldots,x_l\}\subset F_K$.\\
Выберем минимальное $n$ с таким свойством: найдутся $f_i\in F_K$; $0\neq m_i\in {\bf Z}$; не кончающиеся символом из $F_K\cap X$ элементы $v_i$, $\hat{v}_i\in S$ такие, что
\begin{eqnarray}\label{1tm10}
w\equiv (v_1f_1\hat{v}_1^{-1})^{m_1}\ldots (v_nf_n\hat{v}_n^{-1})^{m_n}\mod{(F_K\cap N)^FM},
\end{eqnarray}
где $v_if_i\hat{v}_i^{-1}\in N$ $(i=1,\ldots,n)$.
Числом с таким свойством будет, например, $l$.
Без потери общности рассуждений мы можем и будем считать, что $v_1$ - элемент максимальной длины среди
элементов $v_1, \hat{v}_1,\ldots,v_n, \hat{v}_n$ и ему нет равных в этом множестве элементов.
Действительно, если $\hat{v}_1$ - элемент максимальной длины, то заменим $v_1f_1\hat{v}_1^{-1}$ на $\hat{v}_1f_1^{-1}v_1^{-1}$;
если $v_1=\hat{v}_1$, то $v_1f_1\hat{v}_1^{-1}\in (F_K\cap N)^FM$ --- в противоречии с минимальностью $n$; если
$v_1=v_i$, $1\neq i$, то заменим $v_if_i\hat{v}_i^{-1}$ на $v_1f_1\hat{v}_1^{-1}\hat{v}_1f_1^{-1}f_i\hat{v}_i^{-1}$; если
$v_1=\hat{v}_i$, $1\neq i$, то заменим $v_if_i\hat{v}_i^{-1}$ на $(v_1f_1\hat{v}_1^{-1}\hat{v}_1f_1^{-1}f_i^{-1}v_i^{-1})^{-1}$.\\
Пусть $q\in (I\cup J)\setminus K$ такое, что либо $v_1$ кончается одним из символов $g_q$, $g_q^{-1}$ либо $v_1$ кончается символом из $A_q$. Из (\ref{1tm10}) получаем
\begin{eqnarray}\label{1tm11}
D_q(w)\equiv \sum_{i=1}^n m_i(D_q(v_i)v_i^{-1} - D_q(\hat{v}_i)\hat{v}_i^{-1})\mod{\mathfrak{N}}.
\end{eqnarray}
Если $u\in S$, то нетрудно видеть, что $D_q(u)u^{-1}$ будет суммой элементов вида $\pm t^{-1}$, $t\in S$, $t$ --- начальный отрезок слова $u$. Поэтому из (\ref{1tm11}) следует, что если $v_1$ кончается символом $g_q$, то
$D_q(w)\equiv m_1v_1^{-1} + \mu\mod{\mathfrak{N}}$; если $v_1$ кончается символом из $A_q$, то
$D_q(w)\equiv -m_1v_1^{-1} + \mu\mod{\mathfrak{N}}$; если $v_1=\tilde{v}_1g_q^{-1}$, то
$D_q(w)\equiv -m_1\tilde{v}_1^{-1} + \mu\mod{\mathfrak{N}}$, где $\mu$ --- сумма элементов вида $\pm t^{-1}$, $t\in S$, $t\neq v_1$ и, в случае $v_1=\tilde{v}_1g_q^{-1}$, $t\neq \tilde{v}_1$.
Тогда $D_q(w)\not\equiv 0 ~\mod{\mathfrak{N}}$ --- в противоречии с (\ref{1tm03}).
\end{proof}

\noindent {\bf Следствие} \cite{Rm7}.
{\sl Пусть $F$ --- свободное произведение нетривиальных групп
$A_i~(i\in I)$ и свободной группы с базой $\{g_j | j\in J\}$,
$N$--- нормальная подгруппа в $F$ такая, что $N\cap A_i=1~(i\in I)$. Тогда
$D_k(v)\equiv ~0\mod{{\bf Z}[F]\cdot (N-1)},~k\in I\cup J$, если и только если $v\in [N,N]$}.

\section{Теорема о свободе для произведений групп с одним определяющим соотношением}

\noindent Пусть $F$ --- свободное произведение нетривиальных групп
$A_i~(i\in I)$, $Y$ --- множество элементов группы $F$,
полученное объединением элементов групп $A_i~(i\in I)$,
$P\subseteq I$, $X$ --- множество элементов группы $F$,
полученное объединением элементов групп $A_i~(i\in P)$, $H$ --- подгруппа группы $F$,
порожденная элементами из $X$, $N$ --- нормальная подгруппа в $F$ такая, что $N\cap A_i=1~(i\in I)$.
Так как $N\cap A_i=1~(i\in I)$, то теорема Куроша о подгруппах показывает, что $N$ свободная группа.\\
Назовем смежный класс $L$ группы $F$ по подгруппе $N$ $\alpha$-классом, если $L\cap H\neq \emptyset$ и $\beta$-классом в противном случае.\\
Назовем длиной $\alpha$-класса $L$ длину самого короткого слова из $L\cap H$. В $\alpha$-классах выберем представителей индукцией по длине класса. Выберем пустое слово в качестве представителя для $N$. Если длина $\alpha$-класса $L$ равна 1, то выберем в $L\cap H$ любое слово длины 1 в качестве представителя этого класса. Пусть в $\alpha$-классах длины,
меньшей $l$, представители уже выбраны, т.е. на этих классах уже определена выбирающая
функция $u\rightarrow \bar{u}$. Пусть $L$ --- произвольный $\alpha$-класс длины $l$.
Возьмем в нем какое-нибудь слово $z_1\ldots z_l$, где либо $z_m\in X$ либо $z_m^{-1}\in X$,
и объявим представителем класса $L$ слово $\overline{z_1\ldots z_{l-1}}z_l$.\\
Назовем длиной $\beta$-класса длину самого короткого слова в нем.
В $\beta$-классах выберем представителей индукцией по длине класса.
Если длина $\beta$-класса равна 1, то выберем в нем любое
слово длины 1 в качестве представителя этого класса. Пусть в $\beta$-классах длины,
меньшей $l$, представители уже выбраны, т.е. на этих классах и всех $\alpha$-классах уже определена выбирающая
функция $u\rightarrow \bar{u}$. Пусть $L$ --- произвольный $\beta$-класс длины $l$.
Возьмем в нем какое-нибудь слово $z_1\ldots z_l$, где либо $z_m\in Y$ либо $z_m^{-1}\in Y$,
и объявим представителем класса $L$ слово $\overline{z_1\ldots z_{l-1}}z_l$.
Ясно, что так построенная система представителей, обозначим ее через $S$, шрайерова.
Будем обозначать подсистему в $S$, состоящую из представителей $\alpha$-классов через $S_\alpha$,
состоящую из представителей $\beta$-классов через $S_\beta$.
Элементы из $S_\alpha$ будем называть $\alpha$-представителями, из $S_\beta$ --- $\beta$-представителями.\\
Пусть $N=N_1 \geqslant \ldots \geqslant N_t \geqslant \ldots $ --- нормальный ряд
с абелевыми факторами без кручения, $[N_i,N_j\,]\leqslant N_{i+j}$, $\varphi$ --- естественный гомоморфизм $N\rightarrow N/N_l$, $G$ и $H$ --- подгруппы группы $\varphi(N)$, $G\geqslant H$, $G_i=G\cap \varphi(N_i)$, $H_i=H\cap G_i$.
Положим $\Delta_i$ --- идеал, порожденный $(G_{i_1} - 1)\cdots (G_{i_t} - 1)$
в ${\bf Z}[G]$, $\Delta_i^\prime$ --- идеал, порожденный
$(H_{i_1} - 1)\cdots (H_{i_t} - 1)$ в ${\bf Z}[H]$, где $i_1 + \cdots + i_t\geqslant i$.\\
Если $G$ --- конечно-порожденная (к.п.) подгруппа группы $\varphi(N)$, то мальцевские базы
в $G$ и $H$ будем выбирать согласованными так, чтобы элементы базы группы $H$ из $H_i\setminus H_{i+1}$ были степенями по модулю $G_{i+1}$ некоторых элементов  базы группы $G$.
Мальцевскую базу группы $G$ будем обозначать $M=\{a_1,\ldots ,a_s\}$ и выбирать $M_1,\,M_2\subseteq M$ так, чтобы было $M=M_1\cup M_2$, $M_1\cap M_2=\emptyset$, элементы базы группы $H$ из $H_i\setminus H_{i+1}$ были степенями по модулю $G_{i+1}$ некоторых элементов из $M_1$
и если $m\in M_1$, $m\in G_i\setminus G_{i+1}$, то некоторая степень $m$ по модулю $G_{i+1}$ --- элемент базы группы $H$ из $H_i\setminus H_{i+1}$.\\
Произведение $(a_1-1)^{\alpha_1}\cdots (a_s-1)^{\alpha_s}$, где $\alpha_j$ --- целые неотрицательные числа, называется мономом от $(a_1-1),\ldots,(a_s-1)$. Положим $\omega(a_j-1)=k\Leftrightarrow a_j\in G_k\setminus G_{k+1}$.
Сумма $\alpha_1\cdot\omega(a_1-1) + \ldots + \alpha_s\cdot\omega(a_s-1)$ называется весом монома. Мономы упорядочиваем по весу. Мономы равного веса упорядочиваем по длине. Если $i<j$, то полагаем $(a_i-1)<(a_j-1)$
(полагаем $(a^{n_i}_i-1)<(a^{n_j}_j-1)$ если $a^{n_i}_i$, $a^{n_j}_j$ элементы базы группы $H$). Мономы равного веса и равной длины упорядочиваем лексикографически (слева направо).\\
Непосредственно проверяется, что
\begin{eqnarray}\label{fm01}
n(a-1)\equiv (a^n-1)\mod {\Delta_{k+1}},~a \in G_k\setminus G_{k+1},
\end{eqnarray}
\begin{eqnarray}\label{fm02}
(a-1)(b-1)=(b-1)(a-1)+ba([a,b]-1),~a,\,b \in G.
\end{eqnarray}
Формула (\ref{fm02}) показывает, что любой элемент из $\Delta_{k+1}$ можно записать в виде
целочисленной линейной комбинации элементов из
\begin{eqnarray}\label{psi_0}
\{(a_{i_1}^{\alpha_1}-1)\cdots (a_{i_t}^{\alpha_t}-1)\mid i_1 \leqslant \cdots \leqslant i_t\},
\end{eqnarray}
где $0\neq \alpha_j\in {\bf Z}$, $\omega(a_{i_1}-1) + \ldots + \omega(a_{i_t}-1)\geqslant k+1$.\\
Индукцией по числу $s$ элементов мальцевской базы группы $G$ нетрудно показать, что для любого ненулевого элемента $u\in {\bf Z}[G]$
найдутся такие $k\in {\bf N}$; $m_1,\ldots,m_j$ --- попарно различные мономы веса $k$; $z_1,\ldots,z_j$ --- ненулевые целые числа, что
$u-(z_1m_1+\ldots+z_jm_j)\in {\Delta_{k+1}}$.
Ясно, что при $s=1$ утверждение справедливо. Остается представить $u$ в виде $u=a_1^t(b_0+(a_1-1)b_1+\ldots +(a_1-1)^zb_z)$,
где $t\in {\bf Z}$; $b_0,\,b_1,\ldots,b_z$ --- элементы из ${\bf Z}[G^\prime]$, $G^\prime=\text{гр }(a_2,\ldots,a_s)$ и применить индукционное предположение.\\
Утверждения 1), 2), 3), 4), сформулированные ниже, известны; здесь приведены эскизы их доказательств.\\
1) Пусть $G$ --- к.п. подгруппа группы $\varphi(N)$, тогда мономы веса $k$ образуют базу ${\bf Z}$-модуля $\Delta_k\mod {\Delta_{k+1}}$.\\
%Требуется доказать линейную независимость над ${\bf Z}$ мономов веса $k$ по модулю $\Delta_{k+1}$.
%Предположим противное.
Пусть $m_1,\ldots,m_j$ --- попарно различные мономы веса $k$; $z_1,\ldots,z_j$ --- ненулевые целые числа; $v\in {\Delta_{k+1}}$ и $u=z_1m_1+\ldots+z_jm_j+v=0$.
Если мальцевская база группы $G$ состоит из одного элемента $a$,
то $u=z(a-1)^k+(a-1)^{k+1}v^\prime$, $0\neq z\in {\bf Z}$, $v^\prime\in {\bf Z}[G]$.
Так как $a-1$ не является делителем нуля, то из $u=0$ следует $z + (a-1)v^\prime=0$. Противоречие.\\
Рассмотрим случай $u=z_1m_1(a_s-1)^{k_1}+\ldots+z_jm_j(a_s-1)^{k_j}+v_1+v_2$,
$m_t$ --- мономы, в которых нет подслова $a_s-1$, $k_l\geqslant 0$,
$z_1,\ldots,z_j$ --- ненулевые целые числа, $v_1$ --- целочисленная линейная комбинация элементов из (\ref{psi_0}), у которых
$i_t<s$, $v_2$ --- целочисленная линейная комбинация элементов из (\ref{psi_0}), у которых
$i_t=s$. Обозначим $k_i=\min \{k_1,\ldots,k_j\}$ и пусть $k_i=0$.
Так как мальцевская база группы $G/(a_s)$
состоит из образов элементов $a_1,\ldots ,a_{s-1}$ относительно гомоморфизма $G\rightarrow G/(a_s)$, то из индуктивных соображений можно считать, что $z_i=0$. Противоречие.
Если $k_i>0$, то из $u=0$ следует, что образ $v_1$ относительно гомоморфизма $G\rightarrow G/(a_s)$ равен нулю. Но тогда $v_1=0$,
$v_2=(a_s-1)v_2^\prime$, $v_2^\prime\in {\Delta_{k+1-l}}$.
Так как $a_s-1$ не является делителем нуля, то из $u=0$ следует
$z_1m_1(a_s-1)^{k_1-1}+\ldots+z_im_i(a_s-1)^{k_i-1}+\ldots+z_jm_j(a_s-1)^{k_j-1}+v_2^\prime=0$.
Продолжая аналогичные рассуждения, найдем $v\in {\Delta_{k+1-lk_i}}$ такой, что
$z_1m_1(a_s-1)^{k_1-k_i}+\ldots+z_im_i+\ldots+z_jm_j(a_s-1)^{k_j-k_i}+v=0$, откуда $z_i=0$. Противоречие.\\
2) Справедлива формула ${\bf Z}[H]\cap \Delta_i=\Delta_i^\prime$.\\
Нетрудно видеть, что $\Delta_i^\prime\subseteq {\bf Z}[H]\cap \Delta_i$. Требуется показать, что если $u\in {\bf Z}[H]\cap \Delta_i$, то $u\in \Delta_i^\prime$.
Ясно, что можно ограничиться случаем $G$ --- к.п. группа.\\
Предположим, что $u \equiv z_1m_1+\ldots+z_jm_j \mod{\Delta_{k+1}^\prime}$,
$m_1,\ldots,m_j$ --- попарно различные мономы веса $k$ из ${\bf Z}[H]$, $z_1,\ldots,z_j$ --- целые ненулевые числа, $k<i$.
Формула (\ref{fm01}) показывает, что
$u \equiv z^\prime_1m^\prime_1+\ldots+z^\prime_jm^\prime_j \mod{\Delta_{k+1}}$, где
$m^\prime_1,\ldots,m^\prime_j$ --- попарно различные мономы веса $k$ из ${\bf Z}[G]$, $z^\prime_1,\ldots,z^\prime_j$ --- целые ненулевые числа. Противоречие.\\
3) Если $u\in \Delta_i\setminus \Delta_{i+1}$, $v\in \Delta_t\setminus \Delta_{t+1}$, то
$uv\in \Delta_{i+t}\setminus \Delta_{i+t+1}$.\\
Предположим $uv\in \Delta_{i+t+1}$.
Мы можем и будем считать, что $G$ --- к.п. группа.
Пусть $u\equiv z_1m_1+\ldots+z_jm_j \mod{\Delta_{i+1}}$, $m_1<\ldots<m_j$ --- мономы веса $i$, $z_1,\ldots,z_j$ --- целые ненулевые числа, $v\equiv z^\prime_1m^\prime_1+\ldots+z^\prime_km^\prime_k \mod{\Delta_{t+1}}$, $m^\prime_1<\ldots<m^\prime_k$ --- мономы веса $t$, $z^\prime_1,\ldots,z^\prime_k$ --- целые ненулевые числа.
Формула (\ref{fm02}) показывает, что $uv\equiv \hat{z}_1\hat{m}_1+\ldots+\hat{z}_s\hat{m}_s+z_jz^\prime_k m \mod{\Delta_{i+t+1}}$,
$\hat{m}_1<\ldots<\hat{m}_s<m$ --- мономы веса $i+t$, $\hat{z}_1,\ldots,\hat{z}_s$ --- целые числа. Противоречие.\\
4) Если $u\in \bigcap\limits_{k \in {\bf N}} \Delta_k$, то $u=0$.\\
Предположим противное. Тогда ввиду 1) можем выбрать
к.п. подгруппу $H$ группы
$G$ такую, что $u\in \Delta_i^\prime\setminus \Delta_{i+1}^\prime$.
Найдется к.п. подгруппа $\widehat{H}>H$ группы $G$ такая, что $u$ будет принадлежать идеалу, порожденному
$(\widehat{H}_{i_1} - 1)\cdots (\widehat{H}_{i_t} - 1)$ в ${\bf Z}[\widehat{H}]$, где $i_1 + \cdots + i_t\geqslant i+1$.
Ввиду 2), получим $u\in \Delta_{i+1}^\prime$. Противоречие.\\
Покажем, что если $G=N$, $H$ --- подгруппа группы $N$, $G_i= N_i$, $H_i=H\cap G_i$, $\Delta_i$ --- идеал, порожденный $(G_{i_1} - 1)\cdots (G_{i_t} - 1)$ в ${\bf Z}[G]$, $\Delta_i^\prime$ --- идеал, порожденный
$(H_{i_1} - 1)\cdots (H_{i_t} - 1)$ в ${\bf Z}[H]$, где $i_1 + \cdots + i_t\geqslant i$,
то
\begin{eqnarray}\label{fm2_pr_gr}
\varphi({\bf Z}[H])\cap \varphi(\Delta_i)=\varphi(\Delta_i^\prime).
\end{eqnarray}
Обозначим $\varphi(G)$ через $\overline{G}$, $\varphi(H)$ --- через $\overline{H}$.
Тогда $\overline{G}$ и $\overline{H}$ --- подгруппы группы $\varphi(N)$, $\overline{G}\geqslant \overline{H}$.
Полагаем $\overline{G}_i=\overline{G}\cap \varphi(N_i)$, $\overline{H}_i=\overline{H}\cap \overline{G}_i$;
$\delta_i$, $\delta_i^\prime$ --- идеалы, порожденные
$(\overline{G}_{i_1} - 1)\cdots (\overline{G}_{i_t} - 1)$ в ${\bf Z}[\overline{G}]$ и
$(\overline{H}_{i_1} - 1)\cdots (\overline{H}_{i_t} - 1)$ в ${\bf Z}[\overline{H}]$ соответственно.\\
Так как $\overline{G}=\varphi(N)$, то $\overline{G}_i=\overline{G}\cap \varphi(N_i)=\varphi(N_i)$.
Следовательно, $\delta_i=\varphi(\Delta_i)$.\\
Если $i\geqslant l$, то $\varphi(N_i)=1$; если $i<l$, то $N_i \geqslant N_l$.\\
Поэтому $\varphi(H)\cap \varphi(N_i)=\varphi(H\cap N_i)$.
Следовательно, $\delta^\prime_i=\varphi(\Delta_i^\prime)$.\\
Формула 2) показывает, что ${\bf Z}[\overline{H}]\cap \delta_i=\delta_i^\prime$. Так как ${\bf Z}[\overline{H}]=\varphi({\bf Z}[H])$, то справедлива формула (\ref{fm2_pr_gr}).\\
Если $G$ --- к.п. подгруппа группы $N/N_l$, $\nu\in {\bf Z}[N]$, $\varphi(\nu)\in \Delta_k\setminus \Delta_{k+1}$, то $\bar{\nu}$ будет обозначать линейную комбинацию мономов веса $k$ такую, что $\overline{\nu}\equiv \varphi(\nu)\mod{\Delta_{k+1}}$.
Пусть $H^\prime\leqslant N/N_l$, $m$ --- моном, $H=H^\prime\cap G$. Будем называть $H$-компонентой монома $m$ моном, получающийся вычеркиванием в $m$ тех $a_j-1$, у которых $a_j\in M_2$.
Будем обозначать через $L(m)$ длину $m$, через $L_H(m)$ длину $H$-компоненты $m$.\\
Если $\bar{\nu}=z_1m_1+\ldots+z_km_k$, $z_1,\ldots,z_k$ --- ненулевые целые числа, $m_1,\ldots,m_k$ --- попарно различные мономы, то полагаем $\psi_G(\nu)=\max_t\, (L(m_t)-L_H(m_t))$.\\
Отметим, что если $c\cdot\varphi(\nu)\notin {\bf Z}[H^\prime]+\Delta_{k+1}$, $c\in {\bf N}$, то $\psi_G(\nu)\neq 0$.
Действительно, если $\psi_G(\nu)= 0$, то $\bar{\nu}$ будет линейной комбинацией мономов веса $k$, совпадающих со своей
$H$-компонентой. Формула (\ref{fm01}) показывает, что тогда найдется $c\in {\bf N}$ такое, что
$c\cdot\bar{\nu}\in {\bf Z}[H^\prime]+\Delta_{k+1}$. Поэтому $c\cdot\varphi(\nu)\in {\bf Z}[H^\prime]+\Delta_{k+1}$.\\
Пусть $G$, $\hat{G}$ --- к.п. подгруппы группы $N/N_l$, $\hat{H}=H^\prime\cap \hat{G}$.
По аналогии с $M$, $\Delta_k$ и $\psi_G$ определяем мальцевскую базу $\hat{M}=\{\hat{a}_1,\ldots ,\hat{a}_{\hat{s}}\}=\hat{M}_1\cup \hat{M}_2$ группы $\hat{G}$, идеал $\hat{\Delta}_k$ в ${\bf Z}[\hat{G}]$ и $\psi_{\hat{G}}$.\\
Функция $\psi_G$ обладает следующими свойствами:\\
а) Пусть $\nu$, $\mu\in {\bf Z}[N]$; $\varphi(\nu)$, $\varphi(\mu)\in {\bf Z}[G]$, тогда
$\psi_G(\nu\mu)=\psi_G(\nu)+\psi_G(\mu)$.\\
Из мономов, входящих в линейную комбинацию  $\bar{\nu}$, на которых $\psi_G$ достигает значения $\psi_G(\nu)$, выберем максимальный моном $M_1$. Из мономов, входящих в линейную комбинацию  $\bar{\mu}$, на которых $\psi_G$ достигает значения $\psi_G(\mu)$, выберем максимальный моном $M_2$.\\
Используя формулу (\ref{fm01}) нетрудно доказать, что произведение мономов веса $i$ и $j$, совпадающих со своими $H$-компонентами, равно по модулю $\Delta_{i+j+1}$ линейной комбинации мономов веса $i+j$, каждый из которых совпадает со своей $H$-компонентой.\\
Поэтому $\psi_G(\nu\mu)=\psi_G(M_1M_2)=\psi_G(M_1)+\psi_G(M_2)=\psi_G(\nu)+\psi_G(\mu)$.\\
б) Для произвольного конечного множества $\nu_1,\ldots ,\nu_p$ элементов кольца ${\bf Z}[N]$ можно выбрать
такую к.п. подгруппу $G$ группы $N/N_l$, что $\varphi(\nu_i)\in {\bf Z}[G]$ и если $\hat{G}$ --- к.п. подгруппа группы $N/N_l$, $\hat{G}\geqslant G$, то $\psi_G(\nu_i)=\psi_{\hat{G}}(\nu_i)$ $(i=1,\ldots ,p)$.\\
Действительно, если $a_j\in M_1\cap \varphi(N_k)$, то $a_j$ принадлежит по модулю $\varphi(N_{k+1})$ группе, порожденной
элементами из $\hat{M}_1\cap \varphi(N_k)$, т.е. $\psi_G(\nu)\geqslant\psi_{\hat{G}}(\nu)$, $\varphi(\nu)\in {\bf Z}[G]$.

\begin{lemma}\label{lm6_2_gr}
Пусть $F$ --- свободное
произведение нетривиальных групп $A_i~(i\in I)$, $P\subset I$, $H$ --- свободное произведение групп $A_i~(i\in P)$,
$N$ --- нормальная подгруппа в $F$ такая, что $N\cap A_i=1~(i\in I)$, $F/N$ --- правоупорядочиваемая
группа, $N=N_1 \geqslant \ldots \geqslant N_t \geqslant \ldots $ --- ряд
нормальных подгрупп группы $F$ с абелевыми факторами без кручения, $[N_i,N_j\,]\leqslant N_{i+j}$, $S=S_\alpha\cup S_\beta$ --- система представителей группы $F$ по подгруппе $N$, $\Delta_k$ --- идеал в ${\bf Z}[N]$, порожденный $(N_{i_1} - 1)\cdots (N_{i_s} - 1)$, $i_1 + \cdots + i_s \geqslant k$.\\
Пусть, далее, $g_j, ~f_p,~h_k$ --- элементы из $S_\alpha$ $(g_s\neq g_t,~f_s\neq f_t,~h_s\neq h_t\text{ при }s\neq t)$,
\begin{eqnarray*}
v = \sum_j g_j \mu_j,~\mu_j\in {\bf Z}[H\cap N];
\end{eqnarray*}
\begin{eqnarray*}
r \equiv \sum_p f_p \nu_p\mod{{\bf Z}[F]\cdot\Delta_t},~\nu_p\in \Delta_{t-1}\setminus \Delta_t;
\end{eqnarray*}
\begin{eqnarray*}
w \equiv \sum_k h_k \lambda_k\mod{{\bf Z}[F]\cdot\Delta_{l-t+1}},~\lambda_k\in \Delta_{l-t}\setminus\Delta_{l-t+1}.
\end{eqnarray*}
Если $v \equiv rw\mod{{\bf Z}[F]\cdot\Delta_l}$,
то найдется $c\in {\bf N}$ такое, что $c\cdot\nu_p\in {\bf Z}[H\cap N]+\Delta_t$.
\end{lemma}
\begin{proof}
Обозначим через $\varphi$ естественный гомоморфизм $F\rightarrow F/N_l$, продолженный по линейности на ${\bf Z}[F]$.
Пусть $H^\prime=\varphi(H\cap N)$, $G$ --- к.п. подгруппа группы $N/N_l$, содержащая образы элементов носителей $\mu_j,~\nu_p,~\lambda_k,~\nu_p^{h_k}$ относительно $\varphi$, $G_i=G\cap \varphi(N_i)$, $\delta_t$ --- идеал, порожденный $(G_{i_1} - 1)\cdots (G_{i_k} - 1)$ в ${\bf Z}[G]$, где $i_1 + \cdots + i_k\geqslant t$. Группу $G$ выбираем так, что
$\varphi(\nu)\in \delta_{t-1}\setminus \delta_t$ и
если $\hat{G}$ --- к.п. подгруппа группы $N/N_l$, $\hat{G}\geqslant G$, то $\psi_G(\nu_p)=\psi_{\hat{G}}(\nu_p)$
и $\psi_G(\nu_p^{h_k})=\psi_{\hat{G}}(\nu_p^{h_k})$.\\
Покажем, что $\psi_G(\nu_p)=\psi_G(\nu_p^{h_k})$.
Пусть, для определенности, $\psi_G(\nu_p)\leqslant \psi_G(\nu_p^{h_k})$.
Обозначим $\textup{гр }(a_1,\ldots ,a_s,\,a_1^{h_k},\ldots ,a_s^{h_k})$ через $\hat{G}$.
Ясно, что если $a_j\in M_1\cap \varphi(N_k)$, то $a_j^{h_k}$ принадлежит по модулю $\varphi(N_{k+1})$ группе, порожденной
элементами из $\hat{M}_1\cap \varphi(N_k)$. Отсюда $\psi_G(\nu_p)\geqslant\psi_{\hat{G}}(\nu_p^{h_k})=\psi_G(\nu_p^{h_k})$, т.е. $\psi_G(\nu_p)=\psi_G(\nu_p^{h_k})$.\\
Предположим, найдутся $\nu_p$ такие, что $c\cdot\nu_p\notin {\bf Z}[H\cap N]+\Delta_t$, $c\in {\bf N}$.
Так как $t\leqslant l$, то $c\cdot\varphi(\nu_p)\notin {\bf Z}[H^\prime]+\varphi(\Delta_t$), $c\in {\bf N}$.
Так как $\delta_t\subseteq\varphi(\Delta_t)$, то $c\cdot\varphi(\nu_p)\notin {\bf Z}[H^\prime]+\delta_t$, $c\in {\bf N}$, поэтому $\psi_G(\nu_p)\neq 0$ для таких $\nu_p$.\\
Пусть $M_\nu>0$ и $M_\lambda\geqslant 0$ --- максимальные значения, принимаемые функцией $\psi_G$ на элементах $\nu_p$ и $\lambda_k$ соответственно. Обозначим $\{f_p\mid\psi_G(\nu_p)=M_\nu\}$ и $\{h_k\mid\psi_G(\lambda_k)=M_\lambda\}$ через $A$ и $B$ соответственно. Так как $F/N$ --- правоупорядочиваемая группа, то найдутся
$f_{p_0}\in A$, $h_{k_0}\in B$ такие, что $f_{p_0}h_{k_0}\not\equiv f_ph_k\mod{N}$
при $(p_0,k_0)\neq (p,k)$, $f_p\in A$, $h_k\in B$.
Из $\psi_G(h_{k_0}^{-1}\nu_{p_0} h_{k_0}\lambda_{k_0})=\psi_G(\nu_{p_0}\lambda_{k_0})=M_\nu+M_\lambda$
и $\psi_G(\mu_i)=0$ получаем противоречие.
\end{proof}

\begin{lemma}\label{lm1_2_gr}
Пусть $F$ --- свободное
произведение нетривиальных групп $A_i~(i\in I)$,
$N$ --- нормальная подгруппа в $F$ допускающая элементарные эндоморфизмы группы $F$ такая, что $F/N$ --- упорядочиваемая
группа, $P\subset I$, $H$ --- свободное произведение групп $A_i~(i\in P)$. Пусть, далее, $S=S_\alpha\cup S_\beta$ --- система представителей группы $F$ по подгруппе $N$,
$u\rightarrow \bar{u}$ --- выбирающая функция, $\delta_1,\ldots,\delta_l$, $\mu_1,\ldots,\mu_k$ --- элементы из $S$,
$\delta_iN <\delta_jN$, $\mu_iN <\mu_jN$ при $i <j$.
Тогда если $\{\overline{\mu_1^{-1}\mu_1},\ldots,\overline{\mu_1^{-1}\mu_k}\}\not\subseteq S_\alpha$, то найдутся
$i_0,\,j_0$ такие, что $\overline{\delta_{i_0}\mu_{j_0}}\in S_\beta$ и $\overline{\delta_{i_0}\mu_{j_0}}\neq\overline{\delta_{i}\mu_{j}}$ при $(i_0,j_0)\neq (i,j)$.
\end{lemma}
\begin{proof}
Обозначим через $H^\prime$ свободное произведение групп $A_i~(i\not\in P)$.\\
Полагаем $A=HN/N$, $B$ --- нормальное замыкание группы $H^\prime N/N$ в $F/N$. Очевидно, $F/N=AB$, $A\cap B=1$.
Так как $\overline{\delta_1\mu_1}\neq\overline{\delta_{i}\mu_{j}}$ при $(1,1)\neq (i,j)$, то будем предполагать, что $\overline{\delta_1\mu_1}\in S_\alpha$.
Обозначим $\overline{\delta_i{\delta_1}^{-1}}N$ и $\overline{{\mu_1}^{-1}\mu_j}N$ через $b_ia_i$ и $\hat{a}_j\hat{b}_j$ соответственно,
где $a_i,\hat{a}_j\in A$, $b_i,\hat{b}_j\in B$.
Обозначим максимальный элемент из $\{b_1,\ldots,b_l\}$ через $x$; максимальный элемент из $\{\hat{b}_1,\ldots,\hat{b}_k\}$ через $z$.
Будем считать, для определенности, что $z>1$.
Пусть $b_{i_0}a_{i_0}$ --- максимальный элемент из $\{b_ia_i|b_i=x\}$, $\hat{a}_{j_0}\hat{b}_{j_0}$ --- максимальный элемент из $\{\hat{a}_j\hat{b}_j|\hat{b}_j=z\}$. Элементы $a_i(\overline{\delta_1\mu_1}N)\hat{a}_j$, обозначим их через $t_{ij}$, лежат в $A$, поэтому из $\overline{\delta_{i_0}\mu_{j_0}}N=b_{i_0}t_{i_0 j_0}\hat{b}_{j_0}$; ${\hat{b}_{j_0}}^{t_{i_0 j_0}}>1$;
$b_{i_0}\geqslant 1$ следует
\begin{eqnarray*}
\overline{\delta_{i_0}\mu_{j_0}}N=b_{i_0}{\hat{b}_{j_0}}^{t_{i_0 j_0}}t_{i_0 j_0}\notin A.
\end{eqnarray*}
Из $\overline{\delta_{i_0}\mu_{j_0}}N = \overline{\delta_i\mu_j}N$ вытекает
$t_{i_0 j_0}=t_{ij}$, $b_{i_0}=b_i$, $\hat{b}_{j_0}=\hat{b}_j$; $b_{i_0} a_{i_0}>b_i a_i$, если $i_0\neq i$; $\hat{a}_{j_0} \hat{b}_{j_0} > \hat{a}_j \hat{b}_j$, если $j_0\neq j$; $\overline{\delta_{i_0}\mu_{j_0}}N = b_{i_0} a_{i_0}(\overline{\delta_1\mu_1}N)\hat{a}_{j_0} \hat{b}_{j_0}>b_i a_i(\overline{\delta_1\mu_1}N)\hat{a}_j \hat{b}_j=\overline{\delta_i\mu_j}N$, если $(i_0,j_0)\neq (i,j)$. Следовательно,
$\overline{\delta_{i_0}\mu_{j_0}}N\neq \overline{\delta_i\mu_j}N$,
если $(i_0,j_0)\neq (i,j)$.
\end{proof}

\begin{lemma}\label{lm2_2_gr}
Пусть $X$ --- свободная группа с базой $\{x_j | j\in J\}$, $X_n$ --- $n$-й член нижнего центрального ряда группы $X$,
$\mathfrak{X}$ --- фундаментальный идеал кольца ${\bf Z}[X]$, $v\in X$. Тогда и только тогда $v\in X_n\setminus X_{n+1}$, когда $D_j(v)\in \mathfrak{X}^{n-1} \, (j\in J)$
и найдется элемент $j_0\in J$ такой, что $D_{j_0}(v)\in \mathfrak{X}^{n-1}\setminus  \mathfrak{X}^n$.
\end{lemma}
\begin{proof}
Известно \cite[с.~557]{Fx}, что базу ${\bf Z}$-модуля $\mathfrak{X}^k/\mathfrak{X}^{k+1}$ образуют элементы
вида $(x_{j_1}-1)\ldots (x_{j_k}-1)+\mathfrak{X}^{k+1}$ и что $v\in X_n$ тогда и только тогда, когда $v-1\in \mathfrak{X}^n$.
Следовательно, утверждение леммы вытекает из равенства\\
$v-1=\sum_{j\in J} (x_j-1)D_j(v)$.
\end{proof}

\begin{lemma}\label{tm0_gr}
Пусть $X$ --- свободная группа с базой $\{x_j | j\in J\}$,
$\mathfrak{X}$ --- фундаментальный идеал кольца ${\bf Z}[X]$, $v\in X$, $K\subseteq J$, $F_K=\textup{гр }(x_j| j\in K)$. Тогда
$v$ удовлетворяет условиям
\begin{eqnarray}\label{0tm02_gr}
D_k(v)\equiv ~0\mod{\mathfrak{X}^n},~k\in J\setminus K;~D_k(v)\in{\bf Z}[F_K]\mod{\mathfrak{X}^n},~k\in  K
\end{eqnarray}
если и только если $v\in \textup{гр }(F_K,\,\gamma_{n+1}X)$.
\end{lemma}
\begin{proof}
Лемма \ref{lm2_2_gr} показывает, что из $v\in \text{гр }(F_K,\,\gamma_{n+1}X)$ следуют сравнения (\ref{0tm02_gr}).
Необходимо доказать обратное.\\
Пусть $\varphi\text{: }{\bf Z}[X]\rightarrow {\bf Z}[X]$ --- эндоморфизм, определяемый отображением
$x_j\rightarrow x_j$ при $j\in K$, $x_j\rightarrow 1$ при $j\in J\setminus K$. Обозначим $\bar{v}=\varphi(v)$.
Ясно, что $\bar{v}\in F_K$ и $D_k(v)\equiv D_k(\bar{v})\mod{\mathfrak{X}^n},~k\in  K$.
Так как $D_k(v\bar{v}^{-1})=D_k(v)\bar{v}^{-1}-D_k(\bar{v})\bar{v}^{-1}$, то
\begin{eqnarray*}
D_k(v\bar{v}^{-1})\equiv ~0\mod{\mathfrak{X}^n},~k\in J,
\end{eqnarray*}
т.е. $v\bar{v}^{-1}\in \gamma_{n+1}X$ (лемма \ref{lm2_2_gr}).
\end{proof}
\noindent Индукцией по длине слова $v$ может быть доказана
\begin{lemma}\label{lm5_2_gr}
Пусть $F$ --- свободное произведение нетривиальных групп
$A_i~(i\in I)$, $\{D_i | i\in I\}$ --- соответствующие этой базе производные Фокса кольца ${\bf Z}[F]$. Пусть, далее, $N$ --- свободная подгруппа группы $F$, $v\in N$,
$\{x_z | z\in J\}$ --- база $N$, $\{\partial_z | z\in J\}$ --- соответствующие этой базе производные Фокса кольца ${\bf Z}[N]$. Тогда $D_i(v)=\sum_{z\in J} D_i(x_z)\partial_z(v)$.
\end{lemma}

\noindent
Пусть $F$ --- свободная группа; $N$, $R$ --- нормальные подгруппы группы $F$, $N\geqslant R$, $F/N$ --- правоупорядочиваемая группа,
\begin{eqnarray*}
N=N_{11} \geqslant \ldots \geqslant N_{1,m_1+1}=N_{21} \geqslant \ldots \geqslant N_{s,m_s+1},
\end{eqnarray*}
где $N_{kl}$ --- l-й член нижнего центрального ряда группы $N_{k1}$, $\sqrt{RN_{11}}=N$.\\
Если построена подгруппа $\sqrt{RN_{k1}}$, то
полагаем $\sqrt{RN_{kl}}$ --- множество всех элементов группы $F$, попадающих в некоторой степени в $\gamma_l (\sqrt{RN_{k1}})R$.
Тогда $\sqrt{RN_{kl}}$ --- нормальная подгруппа группы $F$,
$[\sqrt{RN_{kp}}\,,\sqrt{RN_{km}}\,]\leqslant \sqrt{RN_{k,p+m}}$ и $F/\sqrt{RN_{kl}}$ --- правоупорядочиваемая группа.\\
Действительно, $\sqrt{RN_{k1}}/(\gamma_l (\sqrt{RN_{k1}})R)$ --- нильпотентная группа и в ней совокупность периодических
элементов есть подгруппа.\\
Формула $[\sqrt{RN_{kp}}\,,\sqrt{RN_{km}}\,]\leqslant \sqrt{RN_{k,p+m}}$ следует
из того, что в нильпотентной группе без кручения $\sqrt{RN_{k1}}/\sqrt{RN_{k,p+m}}$ если
$x^t y^n=y^n x^t\,(t,\,n\neq 0)$, то $x y=y x$.\\
Утверждение $F/\sqrt{RN_{kl}}$ --- правоупорядочиваемая группа следует из того, что абелевы группы без кручения --- упорядочиваемые группы и расширение правоупорядочиваемой группы посредством правоупорядочиваемой группы будет правоупорядочиваемой
группой.\\
Отметим, что если $F$ --- свободное произведение групп $A_i~(i=1,\ldots,n)$ $(n>1)$,
$1\neq N$ --- нормальная подгруппа в $F$ такая, что $N\cap A_i=1~(i=1,\ldots,n)$,
$F/N$ --- группа без кручения,
то $N$ свободная группа бесконечного ранга.

\begin{proposition}\label{tm4_gr}
Пусть $F$ --- свободное
произведение нетривиальных групп $A_i~(i=1,\ldots,n)$ $(n>2)$,
$1\neq N$ --- нормальная, допускающая элементарные эндоморфизмы группы $F$ подгруппа в $F$ такая, что $N\cap A_i=1~(i=1,\ldots,n)$,
$F/N$ --- группа без кручения. Пусть, далее,
\begin{eqnarray}\label{tm4_0_gr}
N=N_{11} > \ldots > N_{1,m_1+1}=N_{21} > \ldots > N_{s,m_s+1},
\end{eqnarray}
где $N_{ij}$ --- j-й член нижнего центрального ряда группы $N_{i1}$, $s \geqslant 1$, $R$ --- нормальная подгруппа группы $F$, $R\leqslant N$, $H$ --- свободное произведение групп $A_i~(i=1,\ldots,n-1)$.
Если $H\cap RN_{1j}\neq H\cap N_{1j}$, то $H\cap RN_{kl}\neq H\cap N_{kl}\, (N_{kl}\leqslant N_{1j})$.
\end{proposition}
\begin{proof}
Отметим, что $H\cap N_{kl}=\gamma_l(H\cap N_{k1})$.\\
Ясно, что $H\cap N_{kl}\supseteq\gamma_l(H\cap N_{k1})$.
Пусть $\phi$ --- эндоморфизм группы $F$ такой, что
$\phi(a)=1$ $(a\in A_n)$, $\phi(a)=a$  $(a\in A_i,\,i\neq n)$; $u\in H\cap N_{kl}$.
Тогда $\phi(N_{kl})=\gamma_l(H\cap N_{k1})$, поэтому $u=\phi(u)\in \gamma_l(H\cap N_{k1})$.
Следовательно, $\gamma_l(H\cap N_{k1})\supseteq H\cap N_{kl}$.\\
Предположим, $H\cap RN_{1j}\neq H\cap N_{1j}$. Тогда $H\cap N\neq 1$, поэтому $H\cap N$ --- свободная группа бесконечного ранга.
Покажем, что
\begin{eqnarray}\label{pr1_gr}
H\cap RN_{1l}> \gamma_l(H\cap N)\,~(N_{1l}\leqslant N_{1j}).
\end{eqnarray}
По условию, $H\cap RN_{1j}> \gamma_j(H\cap N)$.
Пусть для всех членов ряда~(\ref{tm4_0_gr}) от $N_{1j}$ до $N_{1l}$ включительно формула~(\ref{pr1_gr}) верна ($l\geqslant j$).
Обозначим через $B$ группу $H\cap N$, через $\mathfrak{X}$ --- фундаментальный идеал кольца ${\bf Z}[B]$.
Пусть $v\in (H\cap RN_{1l})\setminus \gamma_l B$.
Выберем в базе группы $B$ элемент $x_1$ такой, что
$\partial_1(v)\notin \mathfrak{X}^{l-1}$, где
$\partial_1$ --- соответствующая элементу $x_1$ производная Фокса кольца ${\bf Z}[B]$.\\
Выберем в базе группы $B$ элемент $x_2\neq x_1$.
Непосредственная проверка показывает, что
$\partial_1([v,x_2])=\partial_1(v)(x_2 - 1) + \partial_1(v)(1-v^{-1}{x_2}^{-1}vx_2)$.\\
Так как $1-v^{-1}{x_2}^{-1}vx_2\in \mathfrak{X}^2$ \cite{Fx}, то $\partial_1([v,x_2])\notin \mathfrak{X}^l$.\\
Следовательно, $[v,x_2]\in (H\cap RN_{1,l+1})\setminus \gamma_{l+1} B$ и
соображения индукции заканчивают доказательство формулы~(\ref{pr1_gr}).\\
Из (\ref{pr1_gr}) следует $H\cap RN_{21}> H\cap N_{21}$.
Остается заметить, что из $H\cap RN_{k1}> H\cap N_{k1}$ следует
$\gamma_{l}(H\cap RN_{k1})> \gamma_{l}(H\cap N_{k1})~(l\in {\bf N})$, т.е.
$H\cap RN_{kl}> H\cap N_{kl}\, (k>1)$.
\end{proof}

\begin{proposition}\label{tm2_gr}
Пусть $F$ --- свободное
произведение нетривиальных групп $A_i~(i=1,\ldots,n)$, $1\neq N$ --- нормальная, допускающая элементарные эндоморфизмы группы $F$ подгруппа в $F$ такая, что $N\cap A_i=1~(i=1,\ldots,n)$, $F/N$ --- правоупорядочиваемая группа. Пусть, далее,
\begin{eqnarray}\label{tm2_0_gr}
N=N_{11} > \ldots > N_{1,m_1+1}=N_{21} > \ldots > N_{s,m_s+1},
\end{eqnarray}
где $N_{kl}=\gamma_l N_{k1}$,
$r\in N_{1t}\backslash N_{1,t+1}$, $R$ --- нормальная подгруппа, порожденная в группе $F$ элементом $r$, $H$ --- свободное произведение групп $A_i~(i=1,\ldots,n-1)$.
Если $H\cap RN_{21}=H\cap N_{21}$, то $H\cap RN_{kl}=H\cap N_{kl}\,(l=1,\ldots, m_k+1,\,k> 1)$.
\end{proposition}
\begin{proof}
\noindent Отметим, что $H\cap N_{kl}=\gamma_l(H\cap N_{k1})$.
Если $H\cap N =1$, то $H\cap RN_{kl} = H\cap N_{kl}=1$, где $N_{kl}$ --- произвольный член ряда {\rm (\ref{tm2_0_gr})}.
Поэтому будем предполагать, что $H\cap N\neq 1$. Нетрудно видеть, что тогда $H\cap N_{kl}\neq 1$.\\
Построим индукцией подгруппы $\sqrt{RN_{kl}}$ группы $F$.
Подгруппа $\sqrt{RN_{21}}$ --- множество всех элементов группы $F$, попадающих в некоторой степени в $RN_{21}$. Пусть построена подгруппа $\sqrt{RN_{k1}}$. Положим $\sqrt{RN_{kl}}$ --- множество всех элементов группы $F$, попадающих в некоторой степени в $\gamma_l (\sqrt{RN_{k1}})R$.
Тогда $\sqrt{RN_{kl}}$ --- нормальная подгруппа группы $F$,
$[\sqrt{RN_{kp}}\,,\sqrt{RN_{km}}\,]\leqslant \sqrt{RN_{k,p+m}}$ и $F/\sqrt{RN_{kl}}$ --- правоупорядочиваемая группа.
Пусть $D_1,\ldots,D_n$ --- производные Фокса
кольца ${\bf Z}[F]$, $\mathfrak{R}_{kl}={\bf Z}[F]\cdot(\sqrt{RN_{kl}}-1)$.\\
Докажем, что $D_n(r)\not\equiv 0\mod \sqrt{RN_{21}}$.
Предположим противное.
Тогда по теореме~\ref{tm1_gr} нашлись бы $v_1,\ldots,v_d$ из $H\cap \sqrt{RN_{21}}$; $f_1,\ldots,f_d$ из $F$
такие, что
\begin{eqnarray*}
r\equiv v_1^{f_1}\cdots v_d^{f_d}\mod [\sqrt{RN_{21}}\,,\sqrt{RN_{21}}\,].
\end{eqnarray*}
Каждый из $v_1,\ldots,v_d$ попадает в некоторой степени в $N_{21}$, поэтому все эти элементы принадлежат $N_{21}$. Следовательно, $r\in N_{1,t+1}$ и мы получаем противоречие. Покажем, что
\begin{eqnarray}\label{tm2_0_0_gr}
H\cap \sqrt{RN_{kl}}=H\cap N_{kl}\,(l=1,\ldots, m_k+1,\,k> 1).
\end{eqnarray}
По условию, $H\cap \sqrt{RN_{21}}=H\cap N_{21}$.
Пусть формула (\ref{tm2_0_0_gr}) справедлива для всех членов ряда (\ref{tm2_0_gr}) от $N_{21}$ до $N_{kl}$
включительно $(N_{21} \geqslant N_{kl}$, $l\leqslant m_k)$,
$S=S_\alpha\cup S_\beta$ --- система представителей группы $F$ по подгруппе $\sqrt{RN_{k1}}$, $\{x_{kz} | z\in J\}$ --- база группы $\sqrt{RN_{k1}}$, $\{\partial_{kz} | z\in J\}$ --- соответствующие этой базе производные Фокса кольца ${\bf Z}[\sqrt{RN_{k1}}\,]$.
Выберем $v\in H\cap (R\gamma_{l+1} (\sqrt{RN_{k1}}))$. Тогда
\begin{eqnarray}\label{tm2_3_gr}
D_m(v)\equiv D_m(r)\cdot  \sum_{p={1}}^{d} k_p\mu_p+\sum_{z\in J} D_m(x_{kz})\partial_{kz}(u)\mod {\mathfrak{R}_{kl}},
\end{eqnarray}
$u\in \gamma_{l+1} (\sqrt{RN_{k1}})$, $k_p\in {\bf Z}[\sqrt{RN_{k1}}\,]$, $\mu_p\in S$ $(\mu_p\neq \mu_t\text{ при }p\neq t)$, $m=1,\ldots,n$. Будем иметь
\begin{eqnarray}\label{tm2_4_gr}
0\equiv D_n(r)\cdot  \sum_{p={1}}^{d} k_p\mu_p+ \sum_{z\in J} D_n(x_{kz})\partial_{kz}(u)\mod {\mathfrak{R}_{kl}}.
\end{eqnarray}
Обозначим через $\mathfrak{X}$ фундаментальный идеал кольца ${\bf Z}[\sqrt{RN_{k1}}\,]$,
через $\Delta_a$ --- идеал, порожденный $\{(\sqrt{RN_{ki_1}} - 1)\cdots (\sqrt{RN_{ki_c}} - 1)|i_1 + \cdots + i_c\geqslant a\}$ в ${\bf Z}[\sqrt{RN_{k1}}\,]$.\\
Из $D_n(r)\not\equiv 0\mod {\mathfrak{R}_{k1}}$ следует, что
$D_n(r) = \sum_{p={1}}^{\hat{d}}  t_p \delta_p$,
где $t_p\in {\bf Z}[\sqrt{RN_{k1}}\,]$, $\delta_p\in S$ $(\delta_p\neq \delta_q\text{ при }p\neq q)$ и
не все элементы из $\{t_1,\ldots,t_{\hat{d}}\}$ принадлежат $\mathfrak{X}$.
Выберем максимальное $l_0$ такое, что $k_p\in \Delta_{l_0}\mod {\mathfrak{R}_{kl}}$, $p=1,\ldots,d$.
Покажем, что $l_0\geqslant l$. Предположим противное.
Пусть $\{\mu_{i_1},\ldots,\mu_{i_a}\}$ --- представители, у которых
$\{k_{i_1},\ldots,k_{i_a}\}\not\subseteq \Delta_{l_0 +1}\mod {\mathfrak{R}_{kl}}$,
$\{\delta_{j_1},\ldots,\delta_{j_b}\}$ --- представители, у которых
$t_{j_1},\ldots,t_{j_b}$ не принадлежат $\mathfrak{X}$.\\
Группа $F/\sqrt{RN_{k1}}$ - правоупорядочиваема, поэтому
среди элементов $\delta_{j_h}\mu_{i_f}$ найдется элемент
$\delta_{j_0}\mu_{i_0}$ такой, что $\delta_{j_0}\mu_{i_0}\sqrt{RN_{k1}}\neq\delta_{j_h}\mu_{i_f}\sqrt{RN_{k1}}$
при $(j_0,i_0)\neq (j_h,i_f)$.
По лемме \ref{lm2_2_gr}, все $\partial_{kz}(u)$ лежат в $\mathfrak{X}^l$,
поэтому, из (\ref{tm2_4_gr}) вытекает существование $c_{i_0}\in {\bf N}$ такого, что
$c_{i_0} k_{i_0}\in \Delta_{l_0 +1}\mod {\mathfrak{R}_{kl}}$. Пришли к противоречию.\\
Следовательно, (\ref{tm2_3_gr}) можно переписать в виде
\begin{eqnarray}\label{tm2_6_gr}
D_m(v)\equiv \sum_{p={1}}^{d_m} \mu_{pm} g_{pm} \mod {\mathfrak{R}_{kl}},~m=1,\ldots,n,
\end{eqnarray}
где $\mu_{pm}\in \Delta_l$, $g_{pm}\in S$.\\
Так как $H\cap N_{k1}\neq 1$,
то $H\cap \sqrt{RN_{k1}}$ --- свободная группа бесконечного ранга.
Если $x$ --- элемент базы группы $H\cap \sqrt{RN_{k1}}$, то найдется $j\in \{1,\ldots,n-1\}$ такое, что
$D_j(x) \not\equiv 0\mod {\mathfrak{R}_{k1}}$.\\
Пусть $x$ --- элемент базы группы $H\cap \sqrt{RN_{k1}}$.
Будем говорить, что $x$ можно поставить в соответствие строку $(M(x),j_x)$, $j_x\in \{1,\ldots,n-1\}$,
$M(x) \in S_\alpha$, если найдутся $g_x\in {\bf N}$, $0\neq \lambda_x\in {\bf Z}$ такие, что
$D_{j_x}(x) \equiv \lambda_x\cdot M(x) + \sum_{p=1}^{g_x} \lambda_p t_p\mod {\mathfrak{R}_{k1}}$, где $\lambda_p\in {\bf Z}$, $t_p \in S_\alpha$, $M(x)\neq t_p$.\\
Нетрудно видеть, что если $z_1,\,z_2$ --- элементы базы группы $H\cap \sqrt{RN_{k1}}$, то
существует последовательность элементарных нильсеновских преобразований,
переводящая $z_1,\,z_2$ в $x_1,\,x_2$, где
в соответствие $x_1$ поставлена строка $(M(x_1),j_{x_1})$ и
эту строку нельзя поставить в соответствие $x_2$.\\
Действительно, если каждому из $z_1,\,z_2$ может быть поставлена в соответствие строка $(M,j)$,
то $D_j(z_i) \equiv \lambda_i\cdot M + \sum_{p=1}^{g_{z_i}} \lambda_{ip} t_{ip}\mod {\mathfrak{R}_{k1}}$, $M\neq t_{ip}$, $\lambda_i\neq 0$, $i=1,\,2$.
Тогда последовательностью элементарных нильсеновских преобразований переводим $z_1,\,z_2$ в такие $x_1,\,x_2$,
что $D_j(x_1) \equiv \lambda\cdot M + \sum_{p=1}^{g_{x_1}} \lambda_p t_p\mod {\mathfrak{R}_{k1}}$;
$\lambda$ --- наибольший общий делитель $\lambda_1,\,\lambda_2$; $M\neq t_p$ и строку $(M,j)$ нельзя поставить в соответствие $x_2$.
Полагаем $(M(x_1),j_{x_1})=(M,j)$.\\
Мы можем и будем считать, что $v$ принадлежит группе, порожденной элементами $x_1,\ldots,x_p$ из базы группы $H\cap \sqrt{RN_{k1}}$,
в соответствие $x_i$ поставлена строка $(M(x_i),j_i)$ и
эту строку нельзя поставить в соответствие $x_a~(i+1\leqslant a\leqslant p)$.
Обозначим через $\partial_1\ldots,\partial_p$ --- производные Фокса кольца ${\bf Z}[H\cap \sqrt{RN_{k1}}]$ соответствующие элементам $x_1,\ldots,x_p$. Покажем, что
\begin{eqnarray}\label{tm2_8_gr}
\partial_z(v)\in  \Delta_l\mod {\mathfrak{R}_{kl}},~ z=1,\ldots,p.
\end{eqnarray}
Предположим противное.  Из $i$ таких, что $\partial_i(v)\not\in \Delta_l\mod {\mathfrak{R}_{kl}}$ выберем наименьшее.
Ввиду леммы \ref{lm5_2_gr}, $D_{j_i}(v)= \sum_{z=1}^p D_{j_i}(x_z)\partial_z(v)$, откуда
\begin{eqnarray}\label{tm2_7_gr}
D_{j_i}(v) \equiv zM(x_i)\partial_i(v)+\sum_{c=1}^g t_c\lambda_c+\sum_{m={1}}^{d} \mu_m g_m\mod {\mathfrak{R}_{kl}},
\end{eqnarray}
где $0\neq z\in {\bf Z}$, $\lambda_c\in {\bf Z}[\sqrt{RN_{k1}}\,]$, $t_c\in S$, $t_c\neq M(x_i)$,
$\mu_m\in \Delta_l$, $g_m\in S$.
Из (\ref{tm2_6_gr}), (\ref{tm2_7_gr}) следует, что $z\partial_i(v)\in \Delta_l\mod {\mathfrak{R}_{kl}}$. Пришли к противоречию.\\
Положим $H_b=H\cap \sqrt{RN_{kb}}$, $b=1,\ldots,l$.
Обозначим через $\Delta_l^\prime$ идеал, порожденный $(H_{i_1} - 1)\cdots (H_{i_c} - 1)$ в ${\bf Z}[H_1]$,  где $i_1 + \cdots + i_c\geqslant l$.
Пусть $\alpha$, $\beta$ --- элементы кольца ${\bf Z}[F]$. Ясно, что $\alpha\equiv \beta\mod {\mathfrak{R}_{kl}}$
тогда и только тогда, когда равны образы элементов $\alpha$, $\beta$ при естественном гомоморфизме $\varphi : {\bf Z}[F]\rightarrow {\bf Z}[F/\sqrt{RN_{kl}}]$.\\
Ввиду (\ref{fm2_pr_gr})
$\varphi({\bf Z}[H_1])\cap \varphi(\Delta_l)=\varphi(\Delta_l^\prime)$, поэтому ${\bf Z}[H_1]\cap \Delta_l\equiv \Delta_l^\prime\mod {\mathfrak{R}_{kl}}$.
Ввиду (\ref{tm2_8_gr}) $\partial_z(v)\in {\bf Z}[H_1]\cap\Delta_l\mod {\mathfrak{R}_{kl}}$,
т.е. $\partial_z(v)\in \Delta_l^\prime\mod {\mathfrak{R}_{kl}},~ z=1,\ldots,p$.\\
По условию, $H\cap \sqrt{RN_{kb}}=H\cap N_{kb}$, поэтому
\begin{eqnarray}\label{tm2_9_gr}
H_b = \gamma_b (H\cap N_{k1}), ~b=1,\ldots,l.
\end{eqnarray}
Пусть $\mathfrak{X}^\prime$ --- фундаментальный идеал кольца ${\bf Z}[H\cap N_{k1}]$.
Из $\gamma_b (H\cap N_{k1})-1\subseteq (\mathfrak{X}^\prime)^b$
и (\ref{tm2_9_gr}) следует, что $\Delta_l^\prime \subseteq (\mathfrak{X}^\prime)^l \mod {\mathfrak{R}_{kl}}$,
т.е. $\partial_z(v)\in (\mathfrak{X}^\prime)^l\mod {\mathfrak{R}_{kl}}$.
Так как $H\cap \sqrt{RN_{kl}}=H\cap N_{kl}=\gamma_l (H\cap N_{k1})$, то $\partial_z(v)\in (\mathfrak{X}^\prime)^l$. Это означает, что $v-1\in (\mathfrak{X}^\prime)^{l+1}$, откуда $v\in \gamma_{l+1} (H\cap N_{k1})$ \cite{Fx},
т.е. $v\in H\cap N_{k,l+1}$.\\
Если $v\in H\cap \sqrt{RN_{k,l+1}}$, то найдется $c\in {\bf N}$
такое, что $v^c\in H\cap (R\gamma_{l+1} (\sqrt{RN_{k1}}))$. Тогда $v^c\in H\cap N_{k,l+1}$, поэтому $v\in H\cap N_{k,l+1}$.
Теперь соображения индукции заканчивают доказательство.
\end{proof}

\begin{proposition}\label{tm5_gr}
Пусть $F$ --- свободное
произведение нетривиальных групп $A_i~(i=1,\ldots,n)$, $1\neq N$ --- нормальная подгруппа в $F$ допускающая элементарные эндоморфизмы группы $F$ такая, что $N\cap A_i=1~(i=1,\ldots,n)$, $F/N$ --- упорядочиваемая группа,
\begin{eqnarray}\label{tm5_0_gr}
N=N_{11} > \ldots > N_{1,m_1+1}=N_{21} > \ldots > N_{s,m_s+1},
\end{eqnarray}
где $N_{kl}=\gamma_l N_{k1}$.
Пусть, далее, $r\in N_{1t}\backslash N_{1,t+1}\,(t\leqslant m_1)$, $R$ --- нормальная подгруппа, порожденная в группе $F$ элементом $r$, $H$ --- свободное произведение групп $A_i~(i=1,\ldots,n-1)$.\\
Если элемент $r$ не сопряжен по модулю $N_{1,t+1}$ ни с каким элементом из $H$, то
$H\cap RN_{1l}=H\cap N_{1l}\,(l=1,\ldots, m_1+1)$.
\end{proposition}
\begin{proof}
\noindent Если $H\cap N=1$, то $H\cap RN_{1l}=H\cap N_{1l}=1\,(l=1,\ldots, m_1+1)$.
Будем предполагать, что $H\cap N\neq 1$.
Отметим, что $H\cap N_{1l}=\gamma_l(H\cap N)$.
Положим $\sqrt{RN_{1l}}$ --- множество всех элементов группы $F$, попадающих в некоторой степени в $RN_{1l}$.
Тогда $\sqrt{RN_{1l}}$ --- нормальная подгруппа группы $F$,
$[\sqrt{RN_{1p}}\,,\sqrt{RN_{1m}}\,]\leqslant \sqrt{RN_{1,p+m}}$ и $F/\sqrt{RN_{1l}}$ --- правоупорядочиваемая группа.
Покажем, что если для всех членов ряда (\ref{tm5_0_gr}) от $N_{1t}$ до $N_{1l}$
включительно $H\cap \sqrt{RN_{1l}}=H\cap N_{1l}$ $(l \geqslant t)$, то
$H\cap \sqrt{RN_{1,{l+1}}}=H\cap N_{1,{l+1}}\,(l\leqslant m_1)$.\\
Пусть $S=S_\alpha\cup S_\beta$ --- система представителей группы $F$ по подгруппе $N$, $D_1,\ldots,D_n$ --- производные Фокса кольца ${\bf Z}[F]$, $\mathfrak{R}_{1l}={\bf Z}[F]\cdot(\sqrt{RN_{1l}}-1)$, $\Delta_l$ --- идеал, порожденный
$\{(\sqrt{RN_{1i_1}} - 1)\cdots (\sqrt{RN_{1i_p}} - 1)|i_1 + \cdots + i_p\geqslant l\}$ в ${\bf Z}[N]$, $\{x_{1z} | z\in J\}$ --- база группы $N$, $\{\partial_{1z} | z\in J\}$ --- соответствующие этой базе производные Фокса кольца ${\bf Z}[N]$.
Выберем $v\in H\cap (RN_{1,l+1})$. Тогда
\begin{eqnarray}\label{tm3_3_gr}
D_m(v)\equiv \sum_{z\in J} D_m(x_{1z})\partial_{1z}(u) + D_m(r)\cdot  \sum_{p={1}}^{d} k_p\mu_p\mod {\mathfrak{R}_{1l}},
\end{eqnarray}
$u\in N_{1,l+1}$, $k_p\in {\bf Z}[N]$, $\mu_p\in S$ $(\mu_p\neq \mu_q\text{ при }p\neq q)$, $m=1,\ldots,n$. Будем иметь
\begin{eqnarray}\label{tm3_4_gr}
\sum_{z\in J} D_n(x_{1z})\partial_{1z}(u) + D_n(r)\cdot  \sum_{p={1}}^{d} k_p\mu_p\equiv 0\mod {\mathfrak{R}_{1l}}.
\end{eqnarray}
Пусть $\mathfrak{X}$ --- фундаментальный идеал кольца ${\bf Z}[N]$.
По лемме \ref{lm2_2_gr}, $\partial_{1z}(u)\in\mathfrak{X}^l$, $\partial_{1z}(v)\in\mathfrak{X}^{l-1}$,
$\partial_{1z}(r)\in\mathfrak{X}^{t-1}$ $(z\in J)$.
Отметим, что $\Delta_a=\mathfrak{X}^a$ при $a\leqslant t$.\\
Если $x$ --- элемент базы группы $N$, то найдется $j$ такое, что
$D_j(x) \not\equiv 0\mod {\mathfrak{R}_{11}}$.
Будем говорить, что $x$ можно поставить в соответствие строку $(M(x),j_x)$, если найдутся $g_x\in {\bf N}$, $0\neq \lambda_x\in {\bf Z}$ такие, что
$D_{j_x}(x) \equiv \lambda_x\cdot M(x) + \sum_{p=1}^{g_x} \lambda_p t_p\mod {\mathfrak{R}_{11}}$, где $\lambda_p\in {\bf Z}$; $M(x),\,t_p \in S$; $M(x)\neq t_p$.\\
Элемент $x$ базы группы $N$ назовем порождающим\\
первого типа, если ему можно поставить в соответствие строку $(M(x),n)$;\\
второго типа, если ему нельзя поставить в соответствие строку $(M(x),n)$,
но можно поставить в соответствие строку $(M(x),j_x)$, $M(x)\in S_\beta$, $j_x\neq n$;\\
третьего типа, если он не является порождающим первого или второго типа.\\
Элементу базы группы $N$ будем ставить в соответствие строку согласно типу этого элемента.
Если $z_1,\ldots,z_b$ --- элементы базы группы $N$, то
существует последовательность элементарных нильсеновских преобразований,
переводящая $z_1,\ldots,z_b$ в $x_1,\ldots,x_b$, где
в соответствие $x_i$ поставлена строка $(M(x_i),j_i)$ и
эту строку нельзя поставить в соответствие $x_a~(i+1\leqslant a\leqslant b)$:\\
а) Применяем элементарные нильсеновские преобразования к порождающим первого типа, переводящие $z_1,\ldots,z_b$ в $x_1^\prime,\ldots,x_b^\prime$, где сначала идут элементы $x_i^\prime$ первого типа,
в соответствие $x_i^\prime$ поставлена строка $(M(x_i^\prime),n)$ и
эту строку нельзя поставить в соответствие $x_a^\prime~(i+1\leqslant a\leqslant b)$, затем порождающие второго и третьего типа.\\
б) Применяем элементарные нильсеновские преобразования к порождающим второго типа, переводящие $x_1^\prime,\ldots,x_b^\prime$ в $x_1^{\prime\prime},\ldots,x_b^{\prime\prime}$, где сначала идут элементы $x_i^{\prime\prime}$ первого и второго типа,
в соответствие $x_i^{\prime\prime}$ поставлена строка $(M(x_i^{\prime\prime}),j_i)$ и
эту строку нельзя поставить в соответствие $x_a^{\prime\prime}~(i+1\leqslant a\leqslant b)$, затем порождающие третьего типа.\\
в) Применяем элементарные нильсеновские преобразования к порождающим третьего типа, переводящие $x_1^{\prime\prime},\ldots,x_b^{\prime\prime}$ в $x_1,\ldots,x_b$, где
в соответствие $x_i$ поставлена строка $(M(x_i),j_i)$ и
эту строку нельзя поставить в соответствие $x_a~(i+1\leqslant a\leqslant b)$.\\
Мы можем и будем считать, что $r$ и $v$ принадлежат группе, порожденной элементами $x_1,\ldots,x_b$ из базы группы $N$,
в соответствие $x_k$ поставлена строка $(M(x_k),j_k)$ и
эту строку нельзя поставить в соответствие $x_a~(k+1\leqslant a\leqslant b)$.
Обозначим через $\partial_1\ldots,\partial_b$ --- производные Фокса кольца ${\bf Z}[N]$ соответствующие элементам $x_1,\ldots,x_b$.\\
Рассмотрим случай $\{\partial_z(v) | z=1,\ldots,b\} \not\subseteq \Delta_l\mod {\mathfrak{R}_{1l}}$.
Выберем минимальное $k$ такое, что $\partial_k(r) \in \Delta_{t-1}\setminus \Delta_t$. Ввиду леммы \ref{lm5_2_gr} $D_{j_k}(r)=\sum_{z=1}^b D_{j_k}(x_z)\partial_z(r)$, откуда
\begin{eqnarray*}
D_{j_k}(r)\equiv zM(x_k)\partial_k(r)+\sum_{p=1}^g \lambda_p t_p\mod{{\bf Z}[F]\cdot\Delta_t},
\end{eqnarray*}
где  $0\neq z\in {\bf Z}$, $\lambda_p\in \Delta_{t-1}$, $t_p\in S$, $t_p\neq M(x_k)$.\\
Т.е. $D_{j_k}(r) = a_1\delta_1+ \cdots +a_q\delta_q$, $a_p \in \Delta_{t-1}$ и $a_p \in \Delta_{t-1}\setminus \Delta_t$ для некоторых $p$; $\delta_p\in S$, $\delta_p\neq \delta_q$ при $p\neq q$. Затем,
выбрав минимальное $k^\prime$ такое, что $\partial_{k^\prime}(v) \in \Delta_{l-1}\setminus \Delta_l$,
аналогичными рассуждениями докажем, что $D_{j_{k^\prime}}(v) = a^\prime_1\delta^\prime_1+ \cdots +a^\prime_{q^\prime}\delta^\prime_{q^\prime}$, $a^\prime_p \in \Delta_{l-1}$ и $a^\prime_p \in \Delta_{l-1}\setminus \Delta_l$ для некоторых $p$, $\delta^\prime_p\in S$, $\delta^\prime_p\neq \delta^\prime_q$ при $p\neq q$.\\
Тогда из (\ref{tm3_3_gr}) следует, что $k_p \in \Delta_{l-t}$, $p\in \{1,\ldots,d\}$ и
$k_p \in \Delta_{l-t}\setminus \Delta_{l-t+1}$ для некоторых $p$.
Пусть $A$ --- подмножество в $\{1,\ldots,d\}$ такое, что $k_p \in \Delta_{l-t}\setminus \Delta_{l-t+1}$ при $p\in A$,
$B$ --- подмножество в $\{1,\ldots,q\}$ такое, что $a_p \in \Delta_{t-1}\setminus \Delta_t$ при $p\in B$.
Так как $D_m(v)$ --- сумма элементов вида $b\theta$, где $\theta$ --- $\alpha$-представитель, $b\in N\cap H$,
то (\ref{tm3_3_gr}) показывает, что не существует $k_0,p_0$  таких, что $\overline{\delta_{k_0}\mu_{p_0}}\in S_\beta$ и $\overline{\delta_{k_0}\mu_{p_0}}\neq\overline{\delta_{k}\mu_{p}}$ при  $(k_0,p_0)\neq (k,p)$, $k_0, k \in A$, $p_0, p \in B$.\\
Из леммы \ref{lm1_2_gr} следует существование $\mu\in S$ такого, что $\overline{\mu^{-1}\mu_p}\in S_\alpha$, $p \in A$, т.е.
\begin{eqnarray}\label{tm3_4_1_1_gr}
\mu^{-1}\cdot  \sum_{p={1}}^{d} k_p\mu_p\equiv b_1\hat{\mu}_1+ \cdots +b_{\hat{d}}\hat{\mu}_{\hat{d}}\mod{{\bf Z}[F]\cdot\Delta_{l-t+1}},
\end{eqnarray}
$\hat{\mu}_p\in S_\alpha$, $\hat{\mu}_p\neq \hat{\mu}_q$ при $p\neq q$, $b_p\in \Delta_{l-t}\setminus \Delta_{l-t+1}$.\\
Покажем, что если $x_z$ --- порождающий первого или второго типа, то
\begin{eqnarray}\label{tm3_4_1_gr}
\partial_z(r^\mu)\in\mathfrak{X}^t,~z\in\{1,\ldots,b\}.
\end{eqnarray}
Предположим противное. Пусть $k$ --- минимальное число, для которого формула (\ref{tm3_4_1_gr}) неверна,
т.е. $\partial_k(r^\mu)\in\mathfrak{X}^{t-1}\setminus \mathfrak{X}^t$.\\
Рассмотрим случай, когда $x_k$ --- порождающий первого типа и ему поставлена в соответствие строка $(M(x_k),n)$.
В этом случае если $x_z$ --- порождающий\\
первого или второго типа и $z<k$, то $\partial_z(r^\mu)\in\mathfrak{X}^t$,\\
если третьего типа, то $D_n(x_z)\equiv 0\mod {\mathfrak{R}_{11}}$, т.е. $D_n(x_z)\partial_z(r^\mu)\in {\bf Z}[F]\cdot\Delta_t$.\\
Тогда из $D_n(r^\mu)=\sum_{z=1}^b D_n(x_z)\partial_z(r^\mu)$ следует, что
\begin{eqnarray*}
D_n(r^\mu)\equiv K\cdot M(x_k)\partial_k(r^\mu)+\sum_{p=1}^g \lambda_p t_p\mod{{\bf Z}[F]\cdot\Delta_t},
\end{eqnarray*}
где  $0\neq K\in {\bf Z}$, $\lambda_p\in \Delta_{t-1}$, $t_p\in S$, $t_p\neq M(x_k)$. Из упорядочиваемости группы $F/N$, следует, что
\begin{eqnarray}\label{tm3_4_2_gr}
D_n(r^\mu)\mu^{-1}\cdot    \sum_{p={1}}^{d} k_p\mu_p\equiv \hat{a}_1\hat{\delta}_1+ \cdots +\hat{a}_c\hat{\delta}_c\mod {\mathfrak{R}_{1l}},
\end{eqnarray}
$\hat{\delta}_p\in S$, $\hat{\delta}_p\neq \hat{\delta}_q$ при $p\neq q$, $\{\hat{a}_1,\ldots,\hat{a}_c\}\subset {\bf Z}[N]$ и $\hat{a}_p \in \Delta_{l-1}\setminus \Delta_l$ для некоторых $p$.
Формула (\ref{tm3_4_2_gr}) противоречит (\ref{tm3_4_gr}).\\
Рассмотрим случай, когда $x_k$ --- порождающий второго типа и ему поставлена в соответствие строка $(M(x_k),j_k)$, $j_k\neq n$, $M(x_k)\in S_\beta$. Тогда
\begin{eqnarray*}
D_{j_k}(r^\mu)\equiv  K\cdot M(x_k)\partial_k(r^\mu)+\sum_{p=1}^g \lambda_p t_p\mod{{\bf Z}[F]\cdot\Delta_t},
\end{eqnarray*}
где  $0\neq K\in {\bf Z}$, $\lambda_p\in \Delta_{t-1}$, $t_p\in S$, $t_p\neq M(x_k)$. Из упорядочиваемости группы $F/N$, следует, что
\begin{eqnarray}\label{tm3_4_3_gr}
D_{j_k}(r^\mu)\mu^{-1}\cdot  \sum_{p={1}}^{d} k_p\mu_p\equiv \hat{a}_1\hat{\delta}_1+ \cdots +\hat{a}_c\hat{\delta}_c\mod {\mathfrak{R}_{1l}},
\end{eqnarray}
$\hat{\delta}_p\in S$, $\hat{\delta}_p\neq \hat{\delta}_q$ при $p\neq q$, $\{\hat{a}_1,\ldots,\hat{a}_c\}\subset {\bf Z}[N]$ и, для некоторых $p$, $\hat{\delta}_p\in S_\beta$, $\hat{a}_p \in \Delta_{l-1}\setminus \Delta_l$.
Формула (\ref{tm3_4_3_gr}) противоречит (\ref{tm3_3_gr}).\\
Полученные противоречия показывают, что формула (\ref{tm3_4_1_gr}) верна.\\
Покажем, что если $x_z$ --- порождающий третьего типа, то найдется натуральное $K_z$ такое, что
$K_z\partial_z(r^\mu)\in{\bf Z}[N\cap H]\mod{\mathfrak{X}^t}$.\\
Предположим противное. Пусть $k$ --- минимальное число, для которого $x_k$ --- порождающий третьего типа и $K\cdot\partial_k(r^\mu)\not\in{\bf Z}[N\cap H]\mod{\mathfrak{X}^t}$, $K\in {\bf N}$.
Элементу $x_k$ поставлена в соответствие строка $(M(x_k),j_k)$, $j_k\neq n$, $M(x_k)\in S_\alpha$;
если $x_z$ --- порождающий третьего типа и $z<k$, то $K_z\cdot\partial_z(r^\mu)\in{\bf Z}[N\cap H]\mod{\mathfrak{X}^t}$;
если $x_z$ --- порождающий первого или второго типа, то ввиду (\ref{tm3_4_1_gr}) $\partial_z(r^\mu)\in\mathfrak{X}^t$.
Тогда
\begin{eqnarray}\label{tm3_4_1_2_gr}
D_{j_k}(r^\mu)\equiv a_0M(x_k)+ a_1\delta_1+ \cdots +a_q\delta_q\mod{{\bf Z}[F]\cdot\Delta_t},
\end{eqnarray}
где $K\cdot a_0\not\in{\bf Z}[N\cap H]\mod{\mathfrak{X}^t}$ для любого натурального $K$;
$M(x_k),\,\delta_p\in S_\alpha$; $\delta_p\neq M(x_k)$; $a_p\in \Delta_{t-1}\setminus \Delta_t$.
Из (\ref{tm3_3_gr}), (\ref{tm3_4_1_1_gr}), (\ref{tm3_4_1_2_gr})следует
\begin{eqnarray*}
D_{j_k}(v)\equiv (a_0M(x_b)+ a_1\delta_1+ \cdots +a_q\delta_q)\cdot  (b_1\hat{\mu}_1+ \cdots +b_{\hat{d}}\hat{\mu}_{\hat{d}})\mod{{\bf Z}[F]\cdot\Delta_l},
\end{eqnarray*}
$b_p\in \Delta_{l-t}\setminus \Delta_{l-t+1}$, $\hat{\mu}_t\in S_\alpha$, $\hat{\mu}_p\neq \hat{\mu}_q$ при $p\neq q$.
Тогда ввиду леммы \ref{lm6_2_gr} найдется натуральное $K$ такое, что $K\cdot a_0\in{\bf Z}[N\cap H]\mod{\mathfrak{X}^t}$.
Противоречие. Следовательно, найдется натуральное $C$ такое, что для всех порождающих третьего типа $x_z$
\begin{eqnarray}\label{tm3_4_1_0_gr}
C\partial_z(r^\mu)\in{\bf Z}[N\cap H]\mod{\mathfrak{X}^t}.
\end{eqnarray}
Из (\ref{tm3_4_1_gr}), (\ref{tm3_4_1_0_gr}) ввиду леммы \ref{tm0_gr} следует, что $r^{C\mu}\in HN_{1,t+1}$.
Так как $N/N_{1,t+1}$ --- группа с однозначным извлечением корней, то $r^\mu\in HN_{1,t+1}$.\\
Противоречие, поэтому
$\{\partial_z(v) | z=1,\ldots,b\} \subseteq \Delta_l\mod {\mathfrak{R}_{1l}}$ и
\begin{eqnarray}\label{tm2_6_gr_2}
D_m(v)\equiv \sum_{p={1}}^{d_m} \mu_{pm} g_{pm} \mod {\mathfrak{R}_{1l}},~m=1,\ldots,n,
\end{eqnarray}
где $\mu_{pm}\in \Delta_l$, $g_{pm}\in S$.\\
Если $y$ --- элемент базы группы $H\cap N$, то найдется $j\in \{1,\ldots,n-1\}$ такое, что
$D_j(y) \not\equiv 0\mod {N}$.
Будем говорить, что $y$ можно поставить в соответствие строку $(M(y),j_y)$, $j_y\in \{1,\ldots,n-1\}$,
$M(y) \in S_\alpha$, если найдутся $g_y\in {\bf N}$, $0\neq \lambda_y\in {\bf Z}$ такие, что
$D_{j_y}(y) \equiv \lambda_y\cdot M(y) + \sum_{p=1}^{g_y} \lambda_p t_p\mod {N}$, где $\lambda_p\in {\bf Z}$, $t_p \in S_\alpha$, $M(y)\neq t_p$.\\
Пусть $\{y_z | z\in J^\prime\}$ --- база группы $H\cap N$.
Мы можем и будем считать, что $v$ принадлежит группе, порожденной элементами $y_1,\ldots,y_c$ из базы группы $H\cap N$, в соответствие $y_k$ поставлена строка $(M(y_k),j_k)$ и
эту строку нельзя поставить в соответствие $y_a~(k+1\leqslant a\leqslant c)$.
Обозначим через $\partial_1^\prime\ldots,\partial_c^\prime$ --- производные Фокса кольца ${\bf Z}[H\cap N]$ соответствующие элементам $y_1,\ldots,y_c$.
Покажем, что
\begin{eqnarray}\label{tm2_8_gr_2}
\partial_z^\prime(v)\in  \Delta_l\mod {\mathfrak{R}_{1l}},~ z=1,\ldots,c.
\end{eqnarray}
Предположим противное. Из $i$ таких, что $\partial_i^\prime(v)\not\in \Delta_l\mod {\mathfrak{R}_{1l}}$ выберем наименьшее.
Ввиду леммы \ref{lm5_2_gr}, $D_{j_i}(v)= \sum_{z=1}^c D_{j_i}(y_z)\partial_z^\prime(v)$, откуда
\begin{eqnarray}\label{tm2_7_gr_2}
D_{j_i}(v) \equiv zM(y_i)\partial_i^\prime(v)+\sum_{p=1}^g t_p\lambda_p+\sum_{m={1}}^{d} \mu_m g_m\mod {\mathfrak{R}_{1l}},
\end{eqnarray}
где $0\neq z\in {\bf Z}$, $\lambda_p\in {\bf Z}[N]$, $t_p\in S$, $t_p\neq M(y_i)$,
$\mu_m\in \Delta_l$, $g_m\in S$.
Из (\ref{tm2_6_gr_2}), (\ref{tm2_7_gr_2}) следует, что $z\partial_i(v)\in \Delta_l\mod {\mathfrak{R}_{kl}}$. Пришли к противоречию.\\
Положим $H_k=H\cap \sqrt{RN_{1k}}$, $k=1,\ldots,l$.
По условию, $H_k=H\cap N_{1k}$, следовательно
\begin{eqnarray}\label{tm2_9_gr0}
H_k = \gamma_k (H\cap N), ~k=1,\ldots,l.
\end{eqnarray}
Обозначим через $\mathfrak{X}^\prime$ фундаментальный идеал кольца ${\bf Z}[H\cap N]$,
через $\Delta_l^\prime$ --- идеал в ${\bf Z}[H_1]$, порожденный $\{(H_{i_1} - 1)\cdots (H_{i_p} - 1)|i_1 + \cdots + i_p\geqslant l\}$.
Известно \cite{Fx}, что
\begin{eqnarray}\label{tm2_11_gr0}
\gamma_k (H\cap N)-1\subseteq (\mathfrak{X}^\prime)^k.
\end{eqnarray}
Пусть $\alpha$, $\beta$ --- элементы кольца ${\bf Z}[F]$. Ясно, что $\alpha\equiv \beta\mod {\mathfrak{R}_{1l}}$
тогда и только тогда, когда образы при естественном гомоморфизме $\varphi : {\bf Z}[F]\rightarrow {\bf Z}[F/\sqrt{RN_{1l}}]$
элементов $\alpha$, $\beta$ равны.\\
Ввиду (\ref{fm2_pr_gr}) $\varphi({\bf Z}[H_1])\cap \varphi(\Delta_l)=\varphi(\Delta_l^\prime)$, поэтому
${\bf Z}[H_1]\cap \Delta_l\equiv \Delta_l^\prime\mod {\mathfrak{R}_{1l}}$.\\
Ввиду (\ref{tm2_8_gr_2}) $\partial^\prime_z(v)\in {\bf Z}[H_1]\cap\Delta_l\mod {\mathfrak{R}_{1l}}$, т.е.
$\partial^\prime_z(v)\in \Delta_l^\prime\mod {\mathfrak{R}_{1l}}$.
Из (\ref{tm2_9_gr0}), (\ref{tm2_11_gr0}) следует, что
$\partial^\prime_z(v)\in (\mathfrak{X}^\prime)^l\mod {\mathfrak{R}_{1l}},~z\in J^\prime$.\\
Так как $H\cap \sqrt{RN_{1l}}=H\cap N_{1l}=\gamma_l (H\cap N)$, то $\partial^\prime_z(v)\in (\mathfrak{X}^\prime)^l,~z\in J^\prime$.
Это означает, что $v-1\in (\mathfrak{X}^\prime)^{l+1}$, следовательно $v\in \gamma_{l+1} (H\cap N)$ \cite{Fx},
т.е. $v\in H\cap N_{1,l+1}$.\\
Если $v\in H\cap \sqrt{RN_{1,l+1}}$, то найдется $c\in {\bf N}$
такое, что $v^c\in H\cap (RN_{1,l+1})$. Поэтому $v\in H\cap N_{1,l+1}$.
Теперь соображения индукции заканчивают доказательство.
\end{proof}
\noindent Из предложений \ref{tm4_gr}, \ref{tm2_gr}, \ref{tm5_gr} вытекает справедливость теоремы~\ref{tm1_pr_gr}.

\section{Теорема о свободе для произведений групп с конечным числом определяющих соотношений}
\noindent Пусть $G$ --- группа. Следующие преобразования матрицы над ${\bf Z}[G]$ назовем элементарными:
\begin{align}
&\text{Перестановка столбцов $i$ и $j$;}\label{gr2_df1}\\
&\text{Перестановка строк $i$ и $j$;}\label{gr2_df2}\\
&\text{Умножение справа }i\text{-й строки на ненулевой элемент из }{\bf Z}[G];\label{gr2_df3}\\
&\text{Прибавление к }j\text{-й строке }i\text{-й строки,}\label{gr2_df4}\\
&\text{умноженной справа на ненулевой элемент из }{\bf Z}[G]\text{, где }i<j.\notag
\end{align}
Если $M=\|a_{kn}\|$ --- матрица над ${\bf Z}[G]$, $\Phi$ --- последовательность
элементарных преобразований матрицы $M$, то $\Phi(M)$ будет обозначать матрицу, полученную из $M$
последовательностью преобразований $\Phi$.
Матрицу $M$ будем называть треугольной ранга $t$, если
$a_{kk}\neq 0\,(k\leqslant t)$, $a_{kn}=0\,(n<k\leqslant t)$, $a_{kn}=0\,(k> t)$.

\begin{lemma}\label{lm4_2_gr_1}
Пусть $F$ --- группа; $N$ --- нормальная подгруппа группы $F$,
$F/N$ --- правоупорядочиваемая разрешимая группа,
$N=N_1 \geqslant \ldots \geqslant N_m\geqslant \ldots$ --- нормальный ряд с абелевыми факторами без кручения, $[N_p\,,N_q\,]\leqslant N_{p+q}$.\\
Пусть, далее, $\phi$ --- естественный гомоморфизм ${\bf Z}[F]\to {\bf Z}[F/N_m]$, $\Delta_t$ --- идеал, порожденный
$\{\phi(N_{i_1} - 1)\cdots \phi(N_{i_c} - 1)|i_1 + \cdots + i_c\geqslant t\}$ в ${\bf Z}[F/N_m]$, $\phi^\prime$ --- естественный гомоморфизм ${\bf Z}[F/N_m]\to {\bf Z}[F/N]$,
$\|a_{kn}\|$ --- $r\times s$ матрица над ${\bf Z}[F/N_m]$, $\|\phi^\prime(a_{kn})\|$ --- треугольная $r\times s$ матрица над ${\bf Z}[F/N]$ ранга $r$, $\psi$ --- функция на ${\bf Z}[F/N_m]$, такая, что $\psi(0)=\infty$; если $\phi^\prime(\alpha)\neq 0$, то $\psi(\alpha)=0$;
если $\alpha\in \Delta_j\setminus \Delta_{j+1}$, то $\psi(\alpha)=j$.
Тогда матрицу $\|a_{kn}\|$ последовательностью элементарных преобразований {\rm (\ref{gr2_df3})}, {\rm (\ref{gr2_df4})}
можно привести к треугольной $r\times s$ матрице $\|b_{kn}\|$ ранга $r$ такой, что $\psi(b_{kk})\leqslant \psi(b_{kn})$.
\end{lemma}
\begin{proof}
\noindent Так как $F/N_m$ --- разрешимая группа без кручения, то в кольце ${\bf Z}[F/N_m]$ выполняется правое условие Оре \cite{Lv}.
Так как $F/N$, $F/N_m$ --- правоупорядочиваемые группы, то в кольцах ${\bf Z}[F/N]$, ${\bf Z}[F/N_m]$ нет делителей нуля.\\
Обозначим через $\Delta_t^\prime$ идеал, порожденный
$\{\phi(N_{i_1} - 1)\cdots \phi(N_{i_c} - 1)|i_1 + \cdots + i_c\geqslant t\}$ в ${\bf Z}[N/N_m]$,
через $S$ --- систему представителей группы $F$ по подгруппе $N$.\\
Пусть $\alpha\equiv\sum_{p=1}^d \mu_p k_p\mod {\Delta_{i+1}}$, $\beta\equiv\sum_{p=1}^q \bar{\mu}_p \bar{k}_p\mod {\Delta_{j+1}}$, где $\mu_p, \bar{\mu}_p\in S$,
$k_p\in \Delta_i^\prime\setminus \Delta_{i+1}^\prime$, $\bar{k}_p\in \Delta_j^\prime\setminus \Delta_{j+1}^\prime$.\\
Группа $F/N$ - правоупорядочиваема, поэтому
среди элементов $\mu_t\bar{\mu}_l$ найдется элемент
$\mu_{t_0}\bar{\mu}_{l_0}$ такой, что $\mu_{t_0}\bar{\mu}_{l_0} N\neq\mu_t\bar{\mu}_l N$
при $(t_0,l_0)\neq (t,l)$.\\
Так как $\mu_{t_0}k_{t_0}\bar{\mu}_{l_0}\bar{k}_{l_0}\in \Delta_{i+j}\setminus \Delta_{i+j+1}$, то $\alpha\beta\in \Delta_{i+j}\setminus \Delta_{i+j+1}$. Отсюда вытекает, что $\psi(\alpha\beta)=\psi(\alpha)+\psi(\beta)$. Ясно, что $\psi(\alpha+\beta)\geqslant\min\{\psi(\alpha),\psi(\beta)\}$, т.е. функция $\psi$ является нормированием на ${\bf Z}[F/N_m]$.\\
По условию, $\phi^\prime(a_{kk})\neq 0$, а если $n<k$, то $\phi^\prime(a_{kn})= 0$,
поэтому $\psi(a_{kk})=0$, а если $n<k$, то $\psi(a_{kn})>0$, $k=1,\ldots,r$.
Обозначим $i$-ю строку матрицы $\|a_{kn}\|$ через $v_i=(a_{i1},\ldots,a_{is})$.
Нужную нам матрицу $\|b_{kn}\|$ будем строить индукцией.\\
Полагаем $b_{1j}=a_{1j}$, $j=1,\ldots,s$, $\bar{v}_1=(b_{11},\ldots,b_{1s})$.
Ясно, что $\psi(b_{11})\leqslant \psi(b_{1j})$, $j=1,\ldots,s$.
Пусть последовательностью элементарных преобразований {\rm (\ref{gr2_df3})}, {\rm (\ref{gr2_df4})}
строк $v_1,\ldots,v_{t-1}$ матрицы $\|a_{kn}\|$
построены строки $\bar{v}_1,\ldots,\bar{v}_{t-1}$ матрицы $\|b_{kn}\|$
такие, что при $n<k$ $b_{kn}= 0$ и $\psi(b_{kk})\leqslant \psi(b_{kn})$, $n=1,\ldots,s$, $a_{t1}=\ldots=a_{tl}=0$.\\
Если $l=t-1$, то полагаем $b_{tj}=a_{tj}$, $j=1,\ldots,s$, $\bar{v}_t=(b_{t1},\ldots,b_{ts})$.
Ясно, что $\psi(b_{tt})\leqslant \psi(b_{tj})$, $j=1,\ldots,s$.\\
Рассмотрим случай $l<t-1$, $a_{t,l+1}\neq 0$, $\phi^\prime(a_{t,l+1})= 0$. Выберем ненулевые $\beta_1,\,\beta_2\in {\bf Z}[F/N_m]$ так, чтобы было $b_{l+1,l+1}\beta_1=-a_{t,l+1}\beta_2$.
Тогда $\bar{v}_{l+1}\beta_1+v_t\beta_2=(c_{t1},\ldots,c_{tr})$, $c_{t1}=\ldots=c_{t,l+1}= 0$.
Из $\psi(b_{l+1,j}\beta_1)\geqslant \psi(b_{l+1,l+1}\beta_1)= \psi(a_{t,l+1}\beta_2) >\psi(\beta_2)$,
$\psi(a_{tt}\beta_2)=\psi(\beta_2)$ следует, что $\psi(c_{tt})=\psi(\beta_2)$,
$\psi(c_{tt})\leqslant \psi(c_{tj})$, $j=1,\ldots,s$. Отметим, что если $j<t$, то $\psi(c_{tt})< \psi(c_{tj})$.
Продолжая аналогичные рассуждения, построим из строки $v_t$ последовательностью элементарных преобразований {\rm (\ref{gr2_df3})}, {\rm (\ref{gr2_df4})} строку $\bar{v}_t=(b_{t1},\ldots,b_{ts})$ такую, что $b_{t1}=\ldots=b_{t,t-1}=0$ и $\psi(b_{tt})\leqslant \psi(b_{tj})$, $j=1,\ldots,s$.
\end{proof}
\noindent В дополнение к лемме \ref{lm4_2_gr_1} отметим, что преобразованиями {\rm (\ref{gr2_df1})}, {\rm (\ref{gr2_df2})},
{\rm (\ref{gr2_df3})}, {\rm (\ref{gr2_df4})} произвольную матрицу $\|a_{kn}\|$ над ${\bf Z}[F/N_m]$ можно привести к треугольной матрице $\|b_{kn}\|$ такой, что $\psi(b_{kk})\leqslant \psi(b_{kn})$.\\
Действительно, элементарными преобразованиями
можно добиться, чтобы на месте $(1,1)$ был элемент $a_{ij}$ такой, что $\psi(a_{ij})=M$,
где $M$ --- минимальное значение, принимаемое функцией $\psi$ на элементах $a_{kn}$,
а затем с его помощью получить нули в первом столбце нашей матрицы.
Продолжая этот процесс применительно к строкам и столбцам у которых номера больше единицы,
мы дойдем в конце концов до треугольной матрицы $\|b_{kn}\|$ такой, что $\psi(b_{kk})\leqslant \psi(b_{kn})$.
\begin{lemma}\label{lm4_2_gr_2}
Пусть $G$ --- разрешимая группа без кручения; $M$  --- матрица над
${\bf Z}[G]$, $\alpha_1,\ldots,\alpha_t,\,\alpha$ --- строки матрицы $M$ ($\alpha_i$ обозначает $i$-ю строку матрицы $M$).
Пусть $\psi$ --- элементарное преобразование матрицы $M$ и если $\psi$ --- преобразование строк матрицы $M$,
то эти строки отличны от $\alpha$; $M^\psi$  --- матрица, полученная из $M$ преобразованием
$\psi$; $\alpha_1^\psi,\ldots,\alpha_t^\psi,\,\alpha^\psi$ --- строки матрицы $M^\psi$.
Если $\alpha$ --- правая линейная комбинация строк $\alpha_1,\ldots,\alpha_t$, то найдется ненулевой элемент $d\in
{\bf Z}[G]$ такой, что $\alpha^\psi d$  --- правая линейная комбинация строк $\alpha_1^\psi,\ldots,\alpha_t^\psi$.
\end{lemma}
\begin{proof}
Нетрудно понять, что если $\psi$ --- одно из преобразований (\ref{gr2_df1}), (\ref{gr2_df2}), (\ref{gr2_df4}), то лемма справедлива.
Пусть $\psi$ --- преобразование (\ref{gr2_df3}),
$M^\psi$ получена из $M$ умножением справа i-й строки матрицы $M$ на ненулевой элемент $a\in {\bf Z}[G]$.
По условию, найдутся элементы $b_1,\ldots,b_t$ из ${\bf Z}[G]$ такие, что
$\alpha_1b_1+\ldots+\alpha_tb_t=\alpha$. Если $b_i= 0$, то берем $d=1$.
Предположим $b_i\neq 0$.  Так как $G$ --- разрешимая группа без кручения, то в кольце ${\bf Z}[G]$ выполняется правое условие Оре \cite{Lv}. Поэтому найдутся ненулевые $c$, $d\in {\bf Z}[G]$ такие, что $ac=b_id$
и мы получим $\alpha_1b_1d+\ldots+\alpha_i ac+\ldots+\alpha_tb_td=\alpha d$.
\end{proof}

\noindent Доказательство теоремы \ref{tm2_pr_gr}.
Положим $\sqrt{RN_{1l}}$ --- множество всех элементов группы $F$, попадающих в некоторой степени в $RN_{1l}$.
Построим индукцией подгруппы $\sqrt{RN_{kl}}$ группы $F$, $k>1$.
Подгруппа $\sqrt{RN_{21}}$ --- множество всех элементов группы $F$, попадающих в некоторой степени в $RN_{21}$. Пусть построена подгруппа $\sqrt{RN_{k1}}$. Положим $\sqrt{RN_{kl}}$ --- множество всех элементов группы $F$, попадающих в некоторой степени в $\gamma_l (\sqrt{RN_{k1}})R$. Тогда $\sqrt{RN_{kl}}$ --- нормальная подгруппа группы $F$, $F/\sqrt{RN_{kl}}$ --- правоупорядочиваемая группа и $[\sqrt{RN_{kl}}\,,\sqrt{RN_{kt}}\,]\leqslant \sqrt{RN_{k,l+t}}$.\\
Естественные гомоморфизмы ${\bf Z}[F]\to {\bf Z}[F/\sqrt{RN_{k,m_k+1}}\,]$ и ${\bf Z}[F/\sqrt{RN_{k,m_k+1}}\,]\to {\bf Z}[F/\sqrt{RN_{k1}}\,]$ обозначим через $\phi_k$ и $\phi^\prime_k$ соответственно.\\
Идеалы, порожденные $\{\phi_k(\sqrt{RN_{ki_1}} - 1)\cdots \phi_k(\sqrt{RN_{ki_c}} - 1)|i_1 + \cdots + i_c\geqslant t\}$ в ${\bf Z}[F/\sqrt{RN_{k,m_k+1}}\,]$ и ${\bf Z}[\sqrt{RN_{k1}}/\sqrt{RN_{k,m_k+1}}\,]$ обозначим через $\Delta_{kt}$ и $\Delta_{kt}^\prime$
соответственно. Определим функцию $\psi_k$ на ${\bf Z}[F/\sqrt{RN_{k,m_k+1}}\,]$ полагая
$\psi_k(0)=\infty$; если $\phi^\prime_k(\alpha)\neq 0$, то $\psi_k(\alpha)=0$;
если $\alpha\in \Delta_{kj}\setminus \Delta_{k,j+1}$, то $\psi_k(\alpha)=j$.
Естественный гомоморфизм ${\bf Z}[F]\to {\bf Z}[F/N]$ обозначим через $\phi_0$.
Пусть $\psi_0$ --- функция на ${\bf Z}[F/N]$ такая, что
$\psi_0(0)=\infty$; если $\alpha\neq 0$, то $\psi_0(\alpha)=0$.\\
При доказательстве леммы \ref{lm4_2_gr_1} было показано, что функция $\psi$ является нормированием на ${\bf Z}[F/N_m]$. Аналогично доказывается, что функция $\psi_k$ является нормированием на ${\bf Z}[F/\sqrt{RN_{k,m_k+1}}\,]$,
т.е. выполняются условия $\psi_k(\alpha\beta)=\psi_k(\alpha)+\psi_k(\beta)$, $\psi_k(\alpha+\beta)\geqslant\min\{\psi_k(\alpha),\psi_k(\beta)\}$, $\alpha,\beta\in {\bf Z}[F/\sqrt{RN_{k,m_k+1}}\,]$.
Ясно, что функция $\psi_0$ является нормированием на ${\bf Z}[F/N]$.\\
Пусть $D_1,\ldots,D_n$ --- производные Фокса кольца ${\bf Z}[F]$.\\
Обозначим $D_j(r_i)$ через $m_{ij}$, через $M$ --- матрицу $\|m_{ij}\|$, через $m_{ij}^{\phi_k}$ --- образ
$m_{ij}$ при $\phi_k$, через $M^{\phi_k}$ --- матрицу $\|m_{ij}^{\phi_k}\|$.
Элементарными преобразованиями приведем матрицу $M^{\phi_0}$ к треугольной матрице $\|m^{(0)}_{ij}\|$.
Обозначим через $M_0$ матрицу $\|m^{(0)}_{ij}\|$, через $t_0$ ранг $M_0$, через $\Phi_0$
последовательность элементарных преобразований матрицы $M$ такую, что $(\Phi_0(M))^{\phi_0}=M_0$.
%Полагаем $I_0$ --- множество номеров $i_1,\ldots,i_{t_0}$ столбцов $m_{i_1},\ldots,m_{i_{t_0}}$ матрицы $M$ таких, что $\Phi_0(m_{i_j})$ --- $j$-й столбец матрицы $\Phi_0(M)$.\\
Ясно, что $\psi_0(m^{(0)}_{ii})\leqslant \psi_0(m^{(0)}_{ij})$.\\
Пусть $k>0$; $\Phi_{k-1}$ --- последовательность элементарных преобразований матрицы $M$ такая, что $(\Phi_{k-1}(M))^{\phi_{k-1}}=M_{k-1}$; $M_{k-1}=\|m^{(k-1)}_{ij}\|$ --- треугольная матрица; $\psi_{k-1}(m^{(k-1)}_{ii})\leqslant \psi_{k-1}(m^{(k-1)}_{ij})$; $t_{k-1}$ --- ранг $M_{k-1}$. Полагаем $M_{k1}=\Phi_{k-1}(M)$.\\
Так как $(M_{k1}^{\phi_k})^{\phi^\prime_k}=M_{k-1}$, то матрица $(M_{k1}^{\phi_k})^{\phi^\prime_k}$ треугольная. Тогда ввиду леммы~\ref{lm4_2_gr_1} найдется такая последовательность $\Phi_{k1}$ элементарных преобразований (\ref{gr2_df3}), (\ref{gr2_df4})
первых $t_{k-1}$ строк матрицы $M_{k1}$, что элементы матрицы $(\Phi_{k1}(M_{k1}))^{\phi_k}=\|b_{ij}\|$ будут
удовлетворять условиям $b_{ij}=0$ при $j< i$ и $\psi_k(b_{ii})\leqslant \psi_k(b_{ij})$, $i=1,\ldots,t_{k-1}$.
Полагаем $M_{k2}=\Phi_{k1}(M_{k1})$.
Последовательностью $\Phi_{k2}$ элементарных преобразований (\ref{gr2_df3}), (\ref{gr2_df4}) с помощью $1\text{-й},\ldots,t_{k-1}\text{-й}$ строк матрицы $M_{k2}$ добьемся того, что в $1\text{-м},\ldots,t_{k-1}\text{-м}$ столбцах матрицы $(\Phi_{k2}(M_{k2}))^{\phi_k}$ под диагональю будут нули.
Полагаем $M_{k3}=\Phi_{k2}(M_{k2})$.
Последовательностью $\Phi_{k3}$ элементарных преобразований строк и столбцов матрицы $M_{k3}$, номера которых больше $t_{k-1}$, добьемся того, что матрица $(\Phi_{k3}(M_{k3}))^{\phi_k}=\|m^{(k)}_{ij}\|$ будет треугольной и $\psi_k(m^{(k)}_{ii})\leqslant \psi_k(m^{(k)}_{ij})$.
Обозначим через $M_k$ матрицу $\|m^{(k)}_{ij}\|$, через $t_k$ ранг $M_k$, через $\Phi_k$
последовательность $\Phi_{k-1}$, $\Phi_{k1}$, $\Phi_{k2}$, $\Phi_{k3}$.
Из индуктивных соображений можно считать, что построены $M_k$, $\Phi_k$, $t_k$, $k=1,\ldots,s$.
Полагаем $I_s$ --- множество $i_1,\ldots,i_{t_s}$ номеров столбцов $m_{i_1},\ldots,m_{i_{t_s}}$ матрицы $M$ таких, что $\Phi_s(m_{i_j})$ --- $j$-й столбец матрицы $\Phi_s(M)$;
$J = \{1,\ldots,n\}\setminus I_s$; $H$ --- группа,
порожденная группами $A_i,~i\in J$.
Если $H\cap N =1$, то $H\cap RN_{kl} = H\cap N_{kl}=1$, где $N_{kl}$ --- произвольный член ряда {\rm (\ref{end_gr_1})}.
Поэтому будем предполагать, что $H\cap N\neq 1$. Нетрудно видеть, что тогда $H\cap N_{kl}\neq 1$,
следовательно $H\cap \sqrt{RN_{k1}}$ и $H\cap \sqrt{RN_{k1}}$ --- свободные группы бесконечного ранга.\\
Так как $t_s\leqslant m$, то $|J|\geqslant n-m$. Пусть $\Phi$ --- последовательность элементарных преобразований матрицы $M$, $\varphi\in\Phi$, $v_1,\ldots,v_n\in {\bf Z}[F]$.
Будем считать, что если $\varphi$ --- преобразование строк матрицы $M$, то $\varphi$ действуют на строку $(v_1,\ldots,v_n)$ тождественно, а если $\varphi$ --- перестановка столбцов $i$ и $j$ матрицы $M$, то $\varphi$ действуют на строку $(v_1,\ldots,v_n)$
перестановкой элементов $v_i$ и $v_j$.\\
Пусть $S_k$ --- система представителей группы $F$ по подгруппе $\sqrt{RN_{k1}}$,
$\mathfrak{X}_k$ --- фундаментальный идеал кольца ${\bf Z}[\sqrt{RN_{k1}}]$,
$\mathfrak{R}_{kl}$ --- идеал, порожденный в ${\bf Z}[F]$ элементами $\{r-1| r\in \sqrt{RN_{kl}}\}$.
Обозначим через $\{x_{kz} \}$
базу группы $\sqrt{RN_{k1}}$, через $\{\partial_{kz} \}$ --- соответствующие этой базе производные Фокса кольца ${\bf Z}[\sqrt{RN_{k1}}\,]$.
Выберем $v\in H\cap (R\gamma_{m_k+1} (\sqrt{RN_{k1}}))$.\\
Найдутся $u\in \gamma_{m_k+1} (\sqrt{RN_{k1}})$, $\beta_1,\ldots,\beta_m\in {\bf Z}[F]$,
такие, что
\begin{eqnarray}\label{tm2_3_gr_end}
D_j(v)\equiv \sum_{i={1}}^{m} D_j(r_i)\beta_i+D_j(u)\mod {\mathfrak{R}_{k,m_k+1}}\,,~j= 1,\ldots,n.
\end{eqnarray}
По лемме \ref{lm2_2_gr}, все $\partial_{kz}(u)$ лежат в $\mathfrak{X}_k^{m_k}$.
Из $D_j(u)= \sum D_j(x_{kz})\partial_{kz}(u)$ следует $D_j(u)\in {\bf Z}[F]\mathfrak{X}_k^{m_k},~j= 1,\ldots,n$. Покажем, что
\begin{eqnarray}\label{tm2_6_gr_end}
\phi_k(D_j(v))\in \Delta_{km_k},~j= 1,\ldots,n.
\end{eqnarray}
Пусть $V=(D_1(v)-D_1(u),\ldots,D_n(v)-D_n(u))$, $\overline{V}=(-D_1(u),\ldots,-D_n(u))$.
Отметим, что элементы строки ${\overline{V}\,}^{\phi_k}$ лежат в $\Delta_{km_k}$.
Ввиду леммы~\ref{lm4_2_gr_2} и (\ref{tm2_3_gr_end}) найдется элемент $d\in {\bf Z}[F]$ такой, что $\phi_k(d)\neq 0$ и строка $(\Phi_k(V d))^{\phi_k}$
линейно выражается через строки  $1,\ldots,t_k$ матрицы $M_k$.\\
Обозначим через $\beta_1,\ldots,\beta_{t_k}$ элементы кольца ${\bf Z}[F]$ такие, что $\phi_k(\beta_i)$ --- коэффициент при
$i$-й строке матрицы $M_k$ в записи $(\Phi_k(V d))^{\phi_k}$ через строки  $1,\ldots,t_k$ матрицы $M_k$.
Первые $t_k$ элементов строки $(\Phi_k(V d))^{\phi_k}$ лежат в $\Delta_{km_k} \phi_k(d)$, так как они соответственно равны первым $t_k$ элементам строки $(\Phi_k(\overline{V} d))^{\phi_k}$.
Так как $\psi_k(m^{(k)}_{ii})\leqslant \psi_k(m^{(k)}_{ij})$, то все элементы строки $(\Phi_k(V d))^{\phi_k}$  лежат в $\Delta_{km_k}\phi_k(d)$, т.е. формула (\ref{tm2_6_gr_end}) справедлива.
(Вначале замечаем, что $m^{(k)}_{11}\phi_k(\beta_1)\in\Delta_{km_k} \phi_k(d)$,
откуда $m^{(k)}_{1j}\phi_k(\beta_1)\in\Delta_{km_k} \phi_k(d),~j= 2,\ldots,n$. Тогда $m^{(k)}_{22}\phi_k(\beta_2)\in\Delta_{km_k} \phi_k(d)$,
откуда $m^{(k)}_{2j}\phi_k(\beta_2)\in\Delta_{km_k} \phi_k(d),~j= 3,\ldots,n$ и т.д.)\\
Ясно, что $H\cap \sqrt{RN_{11}} = H\cap N_{11}$. Предположим, для всех $N_{ij}$ от $N_{11}$ до $N_{kl}$ включительно $H\cap \sqrt{RN_{ij}} = H\cap N_{ij}$, $l\leqslant m_k$.
Возьмем $v\in H\cap (R\gamma_{l+1} (\sqrt{RN_{k1}}))$.\\
Рассмотрим случай $l=1$.
Будем иметь $H\cap \sqrt{RN_{k1}} = H\cap N_{k1}$.\\
Найдутся $u\in \gamma_2 (\sqrt{RN_{k1}})$, $\beta_1,\ldots,\beta_m\in {\bf Z}[F]$, такие, что
\begin{eqnarray*}
D_j(v)\equiv \sum_{i={1}}^{m} D_j(r_i)\beta_i+D_j(u)\mod{\mathfrak{R}_{k1}}\,,~j= 1,\ldots,n.
\end{eqnarray*}
Пусть $V=(D_1(v)-D_1(u),\ldots,D_n(v)-D_n(u))$.
Ввиду леммы~\ref{lm4_2_gr_2} найдется элемент $d\in {\bf Z}[F]$ такой, что $\phi_{k-1}(d)\neq 0$ и строка
$(\Phi_{k-1}(V d))^{\phi_{k-1}}$
линейно выражается через строки треугольной матрицы $M_{k-1}$ ранга $t_{k-1}$.
Так как $u\in \gamma_2 (\sqrt{RN_{k1}})$, то $\phi_{k-1}((D_j(u))=0$ $(j=1,\ldots,n)$,
поэтому в строке $(\Phi_{k-1}(V d))^{\phi_{k-1}}$ первые $t_{k-1}$ элементов --- нули, следовательно,
строка $(\Phi_{k-1}(V d))^{\phi_{k-1}}$ --- нулевая.\\
Тогда $D_j(v)\equiv 0 \mod{N_{k1}}$ $(j=1,\ldots,n)$,
откуда $v\in N_{k2}$ \cite{Rm7}.\\
Если $v\in H\cap \sqrt{RN_{k2}}$, то найдется $c\in {\bf N}$
такое, что $v^c\in H\cap (R\gamma_2 (\sqrt{RN_{k1}}))$. Тогда, как было показано, $v^c\in H\cap N_{k2}$, поэтому $v\in H\cap N_{k2}$.\\
Рассмотрим случай $l>1$. Обозначим через $\{b_z \}$ базу группы $H\cap \sqrt{RN_{k1}}$, через $\{\partial_z \}$ --- производные Фокса, соответствующие этой базе.
Если $x$ --- элемент базы $\{b_z \}$, то найдутся $j_x\in J$, $g_x\in {\bf N}$ такие, что
$D_{j_x}(x) \equiv \sum_{p=1}^{g_x} \lambda_p t_p\mod{\sqrt{RN_{k1}}}$, где $0\neq \lambda_p\in {\bf Z}$, $t_p \in S_\alpha$.
Элементу $x$ поставим в соответствие строку $(M(x),j_x)$, где $M(x)$ --- произвольный элемент из $\{t_1,\ldots,t_{g_x}\}$.
Нетрудно видеть, что если $z_1,\ldots,z_p$ --- элементы базы $\{b_z \}$, то
существует последовательность элементарных нильсеновских преобразований,
переводящая $z_1,\ldots,z_p$ в $x_1,\ldots,x_p$, где
в соответствие $x_i$ поставлена строка $(M(x_i),j_i)$ и
эту строку нельзя поставить в соответствие $x_t~(i+1\leqslant t\leqslant p)$.\\
Мы можем и будем считать, что $v$ принадлежит группе, порожденной элементами $x_1,x_2,\ldots,x_p$ из $\{b_z \}$,
в соответствие $x_i$ поставлена строка $(M(x_i),j_i)$ и
эту строку нельзя поставить в соответствие $x_t~(i+1\leqslant t\leqslant p)$.
Покажем, что
\begin{eqnarray}\label{tm2_8_gr_end}
\phi_k(\partial_z(v))\in \Delta_{kl}^\prime,~ z\in {\bf N}.
\end{eqnarray}
Предположим, формула (\ref{tm2_8_gr_end}) неверна.\\
Из $v\in H\cap N_{kl}$ следует $v\in \gamma_l (H\cap N_{k1})$, поэтому
$\partial_z(v)\in \mathfrak{X}^{l-1},~z=1,2,\ldots,p$.
Выберем наименьшее $i$ такое, что $\phi_k(\partial_i(v))\in \Delta_{k,l-1}^\prime\setminus \Delta_{kl}^\prime$.
Выберем $x_t\in \{x_1,x_2,\ldots,x_p\}$ $(t\neq i)$ и обозначим через $\bar{v}$ элемент $[v,x_t]\in H\cap (R\gamma_{l+2} (\sqrt{RN_{k1}}))$.\\
Тогда $\partial_i(\bar{v})=-\partial_i(v)v^{-1}x_t^{-1}vx_t+\partial_i(v)x_t$. Так как $v^{-1}x_t^{-1}vx_t-1\in \mathfrak{X}^{l+1}$, то
\begin{eqnarray}\label{tm2_8_gr_end0}
\phi_k(\partial_i(\bar{v}))\equiv\phi_k(\partial_i(v)(x_t-1))\mod{\Delta_{k,l+1}^\prime}.
\end{eqnarray}
Обозначим через $\Delta_2$  --- идеал, порожденный $\{(\sqrt{RN_{ki_1}} - 1)\cdots (\sqrt{RN_{ki_c}} - 1)|i_1 + \cdots + i_c\geqslant 2\}$ в ${\bf Z}[\sqrt{RN_{k1}}\,]$, через $\Delta_2^\prime$  --- идеал, порожденный $\{(H\cap \sqrt{RN_{ki_1}} - 1)\cdots (H\cap \sqrt{RN_{ki_c}} - 1)|i_1 + \cdots + i_c\geqslant 2\}$ в $({\bf Z}[H\cap \sqrt{RN_{k1}}\,])$, через $\varphi$ --- гомоморфизм $\phi_k$. Ввиду (\ref{fm2_pr_gr}) $\varphi({\bf Z}[H\cap \sqrt{RN_{k1}}\,])\cap \varphi(\Delta_2)=\varphi(\Delta_2^\prime)$, поэтому
$\phi_k({\bf Z}[H\cap \sqrt{RN_{k1}}\,])\cap\Delta_{k2}^\prime$ --- идеал, порожденный $\{\phi_k(H\cap \sqrt{RN_{ki_1}} - 1)\cdots \phi_k(H\cap \sqrt{RN_{ki_c}} - 1)|i_1 + \cdots + i_c\geqslant 2\}$ в
$\phi_k({\bf Z}[H\cap \sqrt{RN_{k1}}\,])$.\\
Так как $H\cap \sqrt{RN_{k1}}=H\cap N_{k1}$, $H\cap \sqrt{RN_{k2}}=H\cap N_{k2}$, то
$\phi_k({\bf Z}[H\cap \sqrt{RN_{k1}}\,])\cap\Delta_{k2}^\prime$ --- идеал в $\phi_k({\bf Z}[H\cap N_{k1}])$, порожденный $\phi_k(H\cap N_{k1} - 1)^2$ и $\phi_k(H\cap N_{k2} - 1)$. Так как $H\cap N_{k2} = \gamma_2 (H\cap N_{k1})$, то $H\cap N_{k2}-1 \subseteq (H\cap N_{k1} - 1)^2$. Следовательно, идеал в $\phi_k({\bf Z}[H\cap N_{k1}])$, порожденный $\phi_k(H\cap N_{k1} - 1)^2$, равен
$\phi_k({\bf Z}[H\cap N_{k1}])\cap\Delta_{k2}^\prime$.\\
Так как $x_t$ - элемент базы группы $H\cap N_{k1}$, то $x_t-1\not\in {\bf Z}[H\cap N_{k1}]\cdot (H\cap N_{k1}-1)^2$.
Поэтому, ввиду $H\cap \sqrt{RN_{k,m_k+1}}\subset \gamma_2 (H\cap N_{k1})$, $\phi_k(x_t-1)$ не принадлежит
идеалу в $\phi_k({\bf Z}[H\cap N_{k1}])$, порожденному $\phi_k(H\cap N_{k1} - 1)^2$.\\
Следовательно, $\phi_k(x_t-1)\not\in\phi_k({\bf Z}[H\cap N_{k1}])\cap\Delta_{k2}^\prime$, т.е. $\phi_k(x_t-1)\in \Delta_{k1}^\prime\setminus \Delta_{k2}^\prime$.
Так как $\phi_k(\partial_i(v))\in \Delta_{k,l-1}^\prime\setminus \Delta_{kl}^\prime$, то
$\phi_k(\partial_i(v))\phi_k(x_t-1)\in \Delta_{kl}^\prime\setminus \Delta_{k,l+1}^\prime$.\\
Отсюда ввиду (\ref{tm2_8_gr_end0}) вытекает $\phi_k(\partial_i(\bar{v}))\in \Delta_{kl}^\prime\setminus \Delta_{k,l+1}^\prime$.\\
Так как $D_{j_i}(\bar{v})= \sum_{z=1}^p D_{j_i}(x_z)\partial_z(\bar{v})$, то
$D_{j_i}(\bar{v}) = zM(x_i)\partial_i(\bar{v})+\sum_{\bar{p}={1}}^g t_{\bar{p}}\lambda_{\bar{p}}$,
где $0\neq z\in {\bf Z}$, $\lambda_{\bar{p}}\in {\bf Z}[\sqrt{RN_{k1}}\,]$, $t_{\bar{p}}\in S_k$, $t_{\bar{p}}\neq M(x_i)$.
Поэтому $\phi_k(D_{j_i}(\bar{v}))\not\in \Delta_{k,l+1}$.\\
Продолжая аналогичные рассуждения, найдем $w\in H\cap (R\gamma_{m_k+1} (\sqrt{RN_{k1}}))$
такой, что $\phi_k(D_{j_i}(w))\not\in \Delta_{km_k}$
--- противоречие с (\ref{tm2_6_gr_end}).
Следовательно, формула (\ref{tm2_8_gr_end}) верна.
Пусть $H_t=H\cap \sqrt{RN_{kt}}$. Обозначим через
$\Delta_l^\prime$ идеал в ${\bf Z}[H_1]$, порожденный $\{(H_{i_1} - 1)\cdots (H_{i_s} - 1)|i_1 + \cdots + i_s\geqslant l\}$,
$\Delta_l$ --- идеал в ${\bf Z}[\sqrt{RN_{k1}}]$, порожденный $\{(\sqrt{RN_{ki_1}} - 1)\cdots (\sqrt{RN_{ki_c}} - 1)|i_1 + \cdots + i_c\geqslant l\}$.\\
Пусть $\alpha$, $\beta$ --- элементы кольца ${\bf Z}[F]$. Ясно, что $\alpha\equiv \beta\mod {\mathfrak{R}_{kl}}$
тогда и только тогда, когда равны образы элементов $\alpha$, $\beta$ при естественном гомоморфизме $\varphi : {\bf Z}[F]\rightarrow {\bf Z}[F/\sqrt{RN_{kl}}]$.
Ввиду (\ref{fm2_pr_gr})
$\varphi({\bf Z}[H_1])\cap \varphi(\Delta_l)=\varphi(\Delta_l^\prime)$, поэтому
${\bf Z}[H_1]\cap \Delta_l\equiv \Delta_l^\prime\mod {\mathfrak{R}_{kl}}$.\\
Из $\sqrt{RN_{k,m_k+1}}\subset \sqrt{RN_{kl}}$ и (\ref{tm2_8_gr_end}) вытекает, что $\partial_z(v)\in {\bf Z}[H_1]\cap \Delta_l\mod {\mathfrak{R}_{kl}}$, следовательно, $\partial_z(v)\in \Delta_l^\prime\mod {\mathfrak{R}_{kl}},~ z=1,2,\ldots,p$.
По условию, $H_t = \gamma_t (H\cap N_{k1}),~t=1,\ldots,l$, откуда $\partial_z(v)\in \mathfrak{X}^l\mod {\mathfrak{R}_{kl}}$.
Так как $H_l-1\subseteq \mathfrak{X}^l$, то $\partial_z(v)\in \mathfrak{X}^l,~z=1,2,\ldots,p$.
Это означает, что $v-1\in \mathfrak{X}^{l+1}$, следовательно $v\in \gamma_{l+1} (H\cap N_{k1})$ \cite{Fx},
т.е. $v\in H\cap N_{k,l+1}$.\\
Если $v\in H\cap \sqrt{RN_{k,l+1}}$, то найдется $c\in {\bf N}$
такое, что $v^c\in H\cap (R\gamma_{l+1} (\sqrt{RN_{k1}}))$. Тогда $v^c\in H\cap N_{k,l+1}$, поэтому $v\in H\cap N_{k,l+1}$.
Теперь соображения индукции заканчивают доказательство.


\begin{thebibliography}{99}

\bibitem{Fx} R.H.Fox, Free differential calculus. I, Ann. of Math., 57, N~3 (1953), 547--560.

\bibitem{Mg} W.Magnus, \"{U}ber diskontinuierliche Grouppen mit einer definierenden Relation
(Der Freiheitssatz), J. reine u. angew Math., 163, (1931), 141--165.

\bibitem{Rm} Н.С.Романовский, Теорема о свободе для групп с одним определяющим соотношением в многообразиях разрешимых и
нильпотентных групп данных ступеней, Матем. сб., 89(131), N~1(9) (1972), 93--99.

\bibitem{Ya} Г.Г.Ябанжи, Теорема о свободе для групп с одним определяющим соотношением в многообразии
$\mathfrak{N}_c\mathfrak{A}$, Сиб. мат. журн., 21, N~2 (1980), 215--222.

\bibitem{Km} Ю.А.Колмаков, Теорема о свободе для групп с одним определяющим соотношением в многообразии
полинильпотентных групп, Сиб. мат. журн., 27, N~4 (1986), 67--83.

\bibitem{Rm1} Н.С.Романовский, Свободные подгруппы в конечно-определенных группах, Алгебра и логика, 16, N~1 (1977), 88-97.

\bibitem{Rm2} Н.С.Романовский, Обобщенная теорема о свободе для про-$p$-групп, Сиб. матем. ж., 27, N~2 (1986), 154-170.

\bibitem{Rm3} C.K.Gupta, N.S.Romanovski, A generalized Freiheitssatz for centre-by-metabelian groups, Bull. Lond. Math. Soc., 24, N~1 (1992), 71-75.

\bibitem{Rm4} C.K.Gupta, N.S.Romanovski, A generalized Freiheitssatz for the variety $AN_2\cap N_2A$, Algebra Colloq., 1, N~3 (1994), 193-200.

\bibitem{Rm5} Г.Г.Ябанжи, О группах, конечно-определенных в многообразиях $AN_2$ и $N_2A$, Алгебра и логика, 20, N~1 (1981), 109-120.

\bibitem{Rm6} Н.С.Романовский, К теореме о свободе для произведений групп, Алгебра и логика, 38, N~3 (1999), 354-367.

\bibitem{Rm7} Н.С.Романовский, О вложениях Шмелькина для абстрактных и проконечных групп, Алгебра и логика, 38, N~5 (1999), 598-612.

\bibitem{KrMr} М.И.Каргаполов, Ю.И.Мерзляков, Основы теории групп, М.: Наука, 1977, 240 с.

\bibitem{Lv} I.Levin, A note of zero divisors in group rings, Proc. Amer. Math. Soc., 31, N~2 (1972), 357-359.

%\bibitem{GpRm} C.K.Gupta, N.S.Romanovski, On torsion in factors of polynilpotent series of a group with a single relation, Int. J. Algebra Comput., 14, N~4 (2004), 513-523.

%\bibitem{Hr} B.Hartley, The residual nilpotence of wreath products, Proc. London Math. Soc., 20, N~3 (1970), 365-392.

\end{thebibliography}
\end{document}